\documentclass[a4paper,12pt,reqno]{amsart}
\usepackage{amsmath}
\usepackage{amssymb}
\usepackage{amsthm}

\usepackage{amscd}
\usepackage{amsfonts}
\usepackage{setspace}
\usepackage{version}

\usepackage{hyperref}

\title[Typical character and zeta sums]{The typical size of character and zeta sums is $o(\sqrt{x})$}
\author{Adam J Harper}
\address{Mathematics Institute, Zeeman Building, University of Warwick, Coventry CV4 7AL, England}
\email{A.Harper@warwick.ac.uk}
\date{11th January 2023}
\thanks{This work started, with the establishment of a weaker version of Theorem \ref{mainthmchar}, while the author was in residence at the Mathematical Sciences Research Institute in Berkeley, California (supported by the National Science Foundation under Grant No. DMS-1440140), during the Spring 2017 semester. It was also partially supported by the Engineering and Physical Sciences Research Council of the United Kingdom [grant EP/V055755/1].}

\numberwithin{equation}{section}
\theoremstyle{plain}

\onehalfspace
\newcommand{\N}{\mathbb{N}}
\newcommand{\R}{\mathbb{R}}
\newcommand{\E}{\mathbb{E}}

\newcommand{\Z}{\mathbb{Z}}

\newcommand{\Echar}{\mathbb{E}^{\text{char}}}
\newcommand{\Econt}{\mathbb{E}^{\text{cont}}}

\newtheorem{thm1}{Theorem}
\newtheorem{thm2}[thm1]{Theorem}

\newtheorem{thm3}[thm1]{Theorem}

\newtheorem{cor1}{Corollary}
\newtheorem{cor2}[cor1]{Corollary}

\newtheorem{approxres1}{Approximation Result}
\newtheorem{lem1}{Lemma}

\newtheorem{prop1}{Proposition}
\newtheorem{prop2}[prop1]{Proposition}

\newtheorem{harman1}{Harmonic Analysis Result}
\newtheorem{multchaos1}{Multiplicative Chaos Result}

\newtheorem{lem2}[lem1]{Lemma}

\newtheorem{multchaos2}[multchaos1]{Multiplicative Chaos Result}

\begin{document}

\maketitle

\begin{abstract}
We prove conjecturally sharp upper bounds for the Dirichlet character moments $\frac{1}{r-1} \sum_{\chi \; \text{mod} \; r} |\sum_{n \leq x} \chi(n)|^{2q}$, where $r$ is a large prime, $1 \leq x \leq r$, and $0 \leq q \leq 1$ is real. In particular, if both $x$ and $r/x$ tend to infinity with $r$ then $\frac{1}{r-1} \sum_{\chi \; \text{mod} \; r} |\sum_{n \leq x} \chi(n)| = o(\sqrt{x})$, and so the sums $\sum_{n \leq x} \chi(n)$ typically exhibit ``better than squareroot cancellation''. We prove analogous better than squareroot bounds for the moments $\frac{1}{T} \int_{0}^{T} |\sum_{n \leq x} n^{it}|^{2q} dt$ of zeta sums; of Dirichlet theta functions $\theta(1,\chi)$; and of the sums $\sum_{n \leq x} h(n) \chi(n)$, where $h(n)$ is any suitably bounded multiplicative function (for example the M\"{o}bius function $\mu(n)$).

The proofs depend on similar better than squareroot cancellation phenomena for low moments of random multiplicative functions. An important ingredient is a reorganisation of the conditioning arguments from the random case, so that one only needs to ``condition'' on a small collection of fairly short prime number sums. The conditioned quantities arising can then be well approximated by twisted second moments, whose behaviour is the same for character and zeta sums as in the random case.  
\end{abstract}

\section{Introduction}
In this paper we are interested in the size of sums such as
$$ \sum_{n \leq x} n^{it} \;\;\;\;\; \text{and} \;\;\;\;\; \sum_{n \leq x} \chi(n) , $$
where $t \in \R$ and $\chi(n)$ is a non-principal Dirichlet character modulo a large prime $r$. These zeta sums and character sums are among the most studied objects in analytic number theory. We would like to show, on the widest possible range of $x$, that we have substantial cancellation amongst the terms in the sums. Furthermore, we would like to understand the extent of the cancellation.

By periodicity, we can confine our study of character sums to the range $x \leq r$. And if $t$ is large and $x \geq t$, then standard Fourier analysis (see Lemma 1.2 of Ivi\'{c}~\cite{ivic}, for example) shows that $\sum_{x < n \leq 2x} n^{it} = \int_{x}^{2x} w^{it} dw + O(1) = \frac{(2x)^{1+it} - x^{1+it}}{1+it} + O(1)$. So again, we can confine our study of zeta sums to the range $x \leq |t|$.

\vspace{12pt}
For character sums, the classical P\'{o}lya--Vinogradov inequality asserts that we always have $|\sum_{n \leq x} \chi(n)| \ll \sqrt{r} \log r$. This is only non-trivial for $x$ larger than about $\sqrt{r} \log r$, whereas the Burgess bound supplies a non-trivial bound $o(x)$ provided $x \geq r^{1/4 + o(1)}$, for characters modulo prime $r$. See e.g. chapter 9.4 of Montgomery and Vaughan~\cite{mv}. Assuming the Generalised Riemann Hypothesis, Montgomery and Vaughan~\cite{mvexp} improved the P\'{o}lya--Vinogradov inequality to $|\sum_{n \leq x} \chi(n)| \ll \sqrt{r} \log\log r$, and then Granville and Soundararajan~\cite{gransoundlcs} showed that $\sum_{n \leq x} \chi(n) = o(x)$ provided $(\log x)/\log\log r \rightarrow \infty$ as $r \rightarrow \infty$. They also showed this range of $x$ would be best possible for the character sum to always be $o(x)$ (for all such $x$ and all non-principal characters).

Turning to average results, for any $x < r$ we have the easy low moment estimates
$$ \frac{1}{r-1} \sum_{\chi \; \text{mod} \; r} |\sum_{n \leq x} \chi(n)|^2 = \lfloor x \rfloor , \;\;\;\;\; \text{and} \;\;\;\;\; \frac{1}{r-1} \sum_{\chi \; \text{mod} \; r} |\sum_{n \leq x} \chi(n)|^{2q} \leq x^q \;\;\; \forall \; 0 \leq q \leq 1 . $$
The first statement here is a trivial consequence of orthogonality of Dirichlet characters, and the second follows from it using H\"older's inequality. Perhaps surprisingly, these straightforward estimates seem to remain essentially the best known upper bounds for low moments of character sums (i.e. for $0 \leq q \leq 1$). Less straightforwardly, we have
$$ \frac{1}{r-2} \sum_{\chi \neq \chi_0 \; \text{mod} \; r} |\sum_{n \leq x} \chi(n)|^{2q} \ll_{q} r^q \;\; \forall q > 0 , $$
which follows from Theorem 1 of Montgomery and Vaughan~\cite{mvchar} (who actually proved the bound with the fixed sum $\sum_{n \leq x} \chi(n)$ replaced by $M(\chi) := \max_{x} |\sum_{n \leq x} \chi(n)|$). See also Cochrane and Zheng~\cite{cochzheng}, Granville and Soundararajan~\cite{gransoundlcs}, and Kerr~\cite{kerr}, for a selection of stronger upper bounds on high moments when $x$ is small. Notice that it is important to exclude the principal character $\chi_0$ in Montgomery and Vaughan's result~\cite{mvchar} when $q$ and $x$ are large (it would give a large contribution $\asymp x^{2q}/r$), but in the first statement its contribution $\lfloor x \rfloor^{2}/(r-1)$ is not overwhelming (compared with $\lfloor x \rfloor$) provided $x \leq 0.99r$, say.

It is natural to ask whether one typically (e.g. for a positive proportion of characters $\chi$ mod $r$) has squareroot behaviour $\sum_{n \leq x} \chi(n) \asymp \sqrt{x}$, and thus whether the low moments are really $\asymp x^q$ or not\footnote{See the MathOverflow post http://mathoverflow.net/questions/129264/short-character-sums-averaged-on-the-character for explicit discussion of this question.}. Combining the Cauchy--Schwarz inequality with the preceding estimates, for any $x \leq 0.99r$ we get
$$ \frac{1}{r-2} \sum_{\chi \neq \chi_0 \; \text{mod} \; r} |\sum_{n \leq x} \chi(n)|^{2q} \geq \frac{(\frac{1}{r-2} \sum_{\chi \neq \chi_0 \; \text{mod} \; r} |\sum_{n \leq x} \chi(n)|^{2})^2}{\frac{1}{r-2} \sum_{\chi \neq \chi_0 \; \text{mod} \; r} |\sum_{n \leq x} \chi(n)|^{2(2-q)}} \gg_{q} \frac{x^2}{r^{2-q}} \;\;\;\;\; \forall \; 0 \leq q \leq 1 . $$
In particular, for any fixed small $\alpha > 0$ and all $\alpha r \leq x \leq 0.99r$, we now find that indeed $\frac{1}{r-2} \sum_{\chi \neq \chi_0 \; \text{mod} \; r} |\sum_{n \leq x} \chi(n)|^{2q} \asymp_{\alpha, q} x^q$ for all $0 \leq q \leq 1$. Since the order of the second moment is the square of the first moment, another standard Cauchy--Schwarz argument (often called the Paley--Zygmund inequality) implies that $|\sum_{n \leq x} \chi(n)| \gg_{\alpha} \sqrt{x}$ for a positive proportion of characters mod $r$, when $\alpha r \leq x \leq 0.99r$. But when $x=o(r)$, the lower bound $\frac{x^2}{r^{2-q}}$ does not match $x^q$, and the typical size of $\sum_{n \leq x} \chi(n)$ remains unclear. (Note that depending on the size of $x$, one could substantially improve the ``simple'' lower bound $\frac{x^2}{r^{2-q}}$ by suitably applying H\"older's inequality rather than the Cauchy--Schwarz inequality, and using the high moment bounds of e.g. \cite{cochzheng, gransoundlcs, kerr} rather than Montgomery and Vaughan's result~\cite{mvchar}. See also section 1.5 of La Bret\`eche, Munsch and Tenenbaum~\cite{bretechemunschten}. But this would still not deliver a matching lower bound $x^q$ in general.)

\vspace{12pt}
For zeta sums $\sum_{n \leq x} n^{it}$, it is not too difficult to show that $|\sum_{n \leq x} n^{it}| \ll \sqrt{t} \log t$ for all large $x \leq t$, although this estimate seems much less celebrated than the analogous P\'{o}lya--Vinogradov inequality for character sums. See chapter 7.6 of Montgomery~\cite{mont}, and the paper of Fujii, Gallagher and Montgomery~\cite{fgmon}. The Vinogradov--Korobov method yields non-trivial estimates $o(x)$ provided $(\log x)/\log^{2/3}|t| \rightarrow \infty$ as $|t| \rightarrow \infty$, see e.g. chapter 6 of Ivi\'{c}~\cite{ivic}. Obtaining bounds on a wide range of $x$ is particularly important when bounding the Riemann zeta function $\zeta(s)$ for $\Re(s)$ close to 1, with consequences for the distribution of primes. When $\Re(s) = 1/2$, obtaining a larger saving on a more limited range of $x$ is important. The state of the art is recent work of Bourgain~\cite{bourgainzeta}. Somewhat surprisingly, the analogues of the precise conditional character sum results of Montgomery and Vaughan~\cite{mvexp} and of Granville and Soundararajan~\cite{gransoundlcs} do not seem to have been worked out explicitly for zeta sums. However, one could deduce various such results by adapting their methods, using the approximate functional equation for $\zeta(s)$ (which implies in particular that $|\sum_{n \leq x} n^{it}| \asymp \sqrt{t}|\sum_{t/(2\pi x) < n \leq t/\pi} \frac{1}{n^{1+it}} + O(\frac{x}{t} + \frac{1}{\sqrt{t}})|$ for all large $x \leq t$) in place of the P\'olya Fourier expansion for character sums.

Regarding average results, a well known mean value theorem of Montgomery and Vaughan (see e.g. Theorem 5.2 of Ivi\'{c}~\cite{ivic}) asserts that for any complex $a_n$, we have
$$ \frac{1}{T} \int_{0}^{T} |\sum_{n \leq x} a_n n^{it}|^2 dt = \sum_{n \leq x} |a_{n}|^2 (1 + O(\frac{n}{T})) . $$
Thus if $T$ is large and $1 \leq x \leq T$ then $\frac{1}{T} \int_{0}^{T} |\sum_{n \leq x} n^{it}|^2 dt \ll x$, and by H\"older's inequality we get $\frac{1}{T} \int_{0}^{T} |\sum_{n \leq x} n^{it}|^{2q} dt \ll x^q$ for all $0 \leq q \leq 1$. Furthermore, if $T/2 \leq t \leq T$, say, then we can use the approximate functional equation to deduce that $|\sum_{n \leq x} n^{it}| \ll \sqrt{T}|\sum_{T/(2\pi x) < n \leq T/\pi} \frac{1}{n^{1+it}}| + \sqrt{T}$. For any $q \in \N$ we have $(\sum_{T/(2\pi x) < n \leq T/\pi} \frac{1}{n^{1+it}})^q = \sum_{n \leq (T/\pi)^q} \frac{c_{x,T,q}(n)}{n^{1+it}}$, for certain coefficients $c_{x,T,q}(n)$ that are bounded by the $q$-fold divisor function $d_{q}(n)$, and so the mean value theorem yields
$$ \frac{2}{T} \int_{T/2}^{T} |\sum_{n \leq x} n^{it}|^{2q} dt \ll_{q} T^q + T^q \sum_{n \leq (T/\pi)^q} d_{q}(n)^2 (\frac{1}{n^2} + \frac{1}{nT}) \ll_{q} T^q \;\;\; \forall \; q \in \N . $$
By H\"older's inequality we then get this bound for all $q > 0$, which is a partial analogue (with $x$ fixed rather than taking an inner maximum over $x$) of Montgomery and Vaughan's moment bound~\cite{mvchar} for character sums. We remark that although the restriction to $T/2 \leq t \leq T$ could be significantly relaxed here, we would {\em not} have such a bound when integrating over all $0 \leq t \leq T$ with $q$ and $x$ large, due to large contributions from small $t$ (analogously to the contribution from the principal character $\chi_0$ that should be excluded from high moments of character sums).

Now if $x \leq cT$ for a suitable small constant $c > 0$, then Montgomery and Vaughan's mean value theorem also immediately implies that $\frac{2}{T} \int_{T/2}^{T} |\sum_{n \leq x} n^{it}|^2 dt = \lfloor x \rfloor + O(x^{2}/T) \gg x$ (and with more work one could show this for all $1 \leq x \leq T$). So, similarly as for character sums, the Cauchy--Schwarz inequality implies that $\frac{2}{T} \int_{T/2}^{T} |\sum_{n \leq x} n^{it}|^{2q} dt \gg_q x^{2}/T^{2-q}$ for all $1 \leq x \leq T$ and $0 \leq q \leq 1$. In particular, if $x \asymp T$ then these moments are $\asymp x^q$, but when $x=o(T)$ the true order of the low moments remains unclear. See also section 1.5 of La Bret\`eche, Munsch and Tenenbaum~\cite{bretechemunschten} for discussion of lower bounds for the moments.

\vspace{12pt}
One way of exploring the behaviour of $\sum_{n \leq x} \chi(n)$ or $\sum_{n \leq x} n^{it}$ is to consider an appropriate random model. Let $(f(p))_{p \; \text{prime}}$ be a sequence of independent random variables, each distributed uniformly on the complex unit circle. Then we define a {\em Steinhaus random multiplicative function} $f(n)$, by setting $f(n) := \prod_{p^{a} || n} f(p)^{a}$ for all $n \in \N$. Steinhaus random multiplicative functions have been used quite extensively to model a randomly chosen Dirichlet character $\chi(n)$ or ``continuous character'' $n^{it}$: see the papers of Granville and Soundararajan~\cite{gransoundlcs} and Lamzouri~\cite{lamzouri2dzeta}, for example. Helson~\cite{helson} conjectured, by a rough analogy with the first moment of the Dirichlet kernel in classical Fourier analysis, that one should have $\E|\sum_{n \leq x} f(n)| = o(\sqrt{x})$ as $x \rightarrow \infty$. This conjecture was somewhat surprising, given the general philosophy of squareroot (and not more than squareroot) cancellation for oscillating number theoretic sums, and various counter-conjectures were made by other authors. See the introduction to \cite{harperrmflowmoments} for discussion and references. However, the author~\cite{harperrmflowmoments} recently proved Helson's conjecture, in fact showing that
\begin{equation}\label{precisehelson}
\E|\sum_{n \leq x} f(n)|^{2q} \asymp \Biggl(\frac{x}{1 + (1-q)\sqrt{\log\log x}} \Biggr)^q
\end{equation}
uniformly for all large $x$ and all $0 \leq q \leq 1$. This raises the question whether one should now expect better than squareroot cancellation for character and zeta sums, on some range of $x$ rather smaller than the conductor, or whether the random multiplicative model simply fails to capture the arithmetic truth in these problems.

\subsection{Statement of results}
Our main results establish, for character and zeta sums, the natural analogues of the upper bound part of \eqref{precisehelson}.

\begin{thm1}\label{mainthmchar}
Let $r$ be a large prime. Then uniformly for any $1 \leq x \leq r$ and any $0 \leq q \leq 1$, we have
$$ \frac{1}{r-1} \sum_{\chi \; \text{mod} \; r} |\sum_{n \leq x} \chi(n)|^{2q} \ll \Biggl(\frac{x}{1 + (1-q)\sqrt{\log\log(10L)}} \Biggr)^q , $$
where $L = L_r := \min\{x,r/x\}$.
\end{thm1}

\begin{thm2}\label{mainthmcont}
Let $T$ be a large real number. Then uniformly for any $1 \leq x \leq T$ and any $0 \leq q \leq 1$, we have
$$ \frac{1}{T} \int_{0}^{T} |\sum_{n \leq x} n^{it}|^{2q} dt \ll \Biggl(\frac{x}{1 + (1-q)\sqrt{\log\log(10L_T)}} \Biggr)^q  , $$
where $L_T := \min\{x,T/x\}$.
\end{thm2}

In particular, for any {\em fixed} $0 < q < 1$ and any $x=x(r)$ such that $x$ and $r/x$ both tend to infinity with $r$, we have $\frac{1}{r-1} \sum_{\chi \; \text{mod} \; r} |\sum_{n \leq x} \chi(n)|^{2q} \ll \frac{x^q}{(\log\log(10L))^{q/2}} = o(x^q)$.

We shall discuss the proofs in detail in section \ref{subsecproofideas}, below. We note here that there is a well known ``symmetry'' in the behaviour of character sums and zeta sums, whereby e.g. $\frac{1}{\sqrt{x}} |\sum_{n \leq x} \chi(n)| \approx \sqrt{\frac{x}{r}} |\sum_{n \leq r/x} \chi(n)|$ (very roughly speaking). See section 10 of Granville and Soundararajan~\cite{gransoundlcs}. This symmetry is sometimes called the ``Fourier flip'', and manifests itself in the (approximate) functional equations of the corresponding $L$-functions, and in the very structured nature of long sums $\sum_{n \leq x} \chi(n) , \sum_{n \leq x} n^{it}$ where $x$ is of the same order as the conductor (i.e. as $r$ or $|t|$, respectively). Given the symmetry between character sums of lengths $x$ and $r/x$, and between zeta sums of lengths $x$ and $|t|/x$, the quantities $\log\log(10L_r)$ and $\log\log(10L_T)$ appearing in Theorems \ref{mainthmchar} and \ref{mainthmcont} are natural substitutes for the $\log\log x$ saving factor in \eqref{precisehelson}.

The shape of the bounds in Theorems \ref{mainthmchar} and \ref{mainthmcont} might initially seem peculiar, and perhaps open to improvement. In most number theoretic settings, if one obtains a saving one expects to save at least a power of a logarithm. But in fact it seems reasonable to conjecture that Theorems \ref{mainthmchar} and \ref{mainthmcont} are sharp (provided in Theorem \ref{mainthmchar} that $x \leq 0.99r$, say, so we are away from the point where periodicity trivially induces substantial extra cancellation\footnote{For non-principal $\chi$, if $0.99r < x \leq r$ we can observe that $\sum_{n \leq x} \chi(n) = - \sum_{x < n \leq r} \chi(n) = - \sum_{1 \leq n < r-x} \chi(r-n) = - \chi(-1) \sum_{1 \leq n < r-x} \chi(n)$, and then apply Theorem \ref{mainthmchar} to these sums instead.}). Note that in the probabilistic setting of \eqref{precisehelson} we already have an order of magnitude result, rather than just an upper bound.

It is possible that the methods leading to Theorems \ref{mainthmchar} and \ref{mainthmcont} would produce matching lower bounds when $x \leq e^{\log^{c}r}$ and $x \leq e^{\log^{c}T}$, say, for a certain small $c > 0$, although the author has not checked all details of this. (See the discussion of the proofs of Theorems \ref{mainthmchar} and \ref{mainthmcont} in section \ref{subsecproofideas}. To produce lower bounds in an analogous way, one would need to compare {\em conditional second and fourth moments}, rather than simply using H\"older's inequality to pass to conditional second moments, with the parameter $\log P$ in the proof being kept comparable to $\log x$ to prevent blow-up from primes larger than $P$ in the conditional fourth moments. This seems to correspond to conditions like $x \leq e^{\log^{c}r}$ and $x \leq e^{\log^{c}T}$.) By ``symmetry'', this should also lead to lower bounds when $x$ is close to $r$ and $T$, respectively. For general $x$ and $q < 1$, the best existing lower bounds for these moments seem to differ from our upper bounds by powers of $\log L$ (see e.g. section 1.5 of La Bret\`eche, Munsch and Tenenbaum~\cite{bretechemunschten}), and it is a challenging and interesting open problem to show that Theorems \ref{mainthmchar} and \ref{mainthmcont} are sharp. As we shall discuss below in the context of Corollary \ref{cortheta}, matching lower bounds for Theorem \ref{mainthmchar} when $x \approx \sqrt{r}$ (say) would have some important consequences.

\vspace{12pt}
Theorem \ref{mainthmchar} implies that for most characters $\chi$ mod $r$, we have $|\sum_{n \leq x} \chi(n)| \ll \frac{\sqrt{x}}{(\log\log(10L))^{1/4}}$. This is ``better than squareroot cancellation'', provided $L$ is large. The following Corollary, which is the character sum version of Corollary 2 of Harper~\cite{harperrmflowmoments} from the random multiplicative case, makes this observation a little more precise. 

\begin{cor1}\label{cordev}
Let $r$ be a large prime. Uniformly for any $1 \leq x \leq r$ and any $\lambda \geq 2$, we have
$$ \frac{1}{r-1} \#\{\chi \; \text{mod} \; r : |\sum_{n \leq x} \chi(n)| \geq \lambda\frac{\sqrt{x}}{(\log\log(10L))^{1/4}}\} \ll \frac{\min\{\log\lambda, \sqrt{\log\log(10L)}\}}{\lambda^2} , $$
where $L = \min\{x,r/x\}$.
\end{cor1}

Corollary \ref{cordev} follows from Theorem \ref{mainthmchar} with $q$ chosen suitably in terms of $\lambda$. See section \ref{secotherproofs}, below, for the full (short) argument. If desired, one could use Theorem \ref{mainthmcont} to deduce an analogous corollary for zeta sums.

\vspace{12pt}
One of the many motivations for considering the character sums $\sum_{n \leq x} \chi(n)$ is their connection with Dirichlet theta functions $\theta(s,\chi)$. Recall that whenever $\Re(s) > 0$, we define
\begin{equation}
\theta(s,\chi) := \left\{ \begin{array}{ll}
     \sum_{n=1}^{\infty} \chi(n) e^{-\pi n^{2} s/r} & \text{if} \; \chi \; \text{is even} , \\
     \sum_{n=1}^{\infty} n \chi(n) e^{-\pi n^{2} s/r} & \text{if} \; \chi \; \text{is odd} .
\end{array} \right. \nonumber
\end{equation}
(These differ by a factor of 2 from the definitions sometimes given, but that will be unimportant for our purposes here.) Theta functions arise in the theory of automorphic forms, and are an important tool in classical proofs of the functional equation for Dirichlet $L$-functions. See e.g. chapter 10.1 of Montgomery and Vaughan~\cite{mv}.

In the last few years, a sequence of papers have investigated the moments of $\theta(1,\chi)$, as $\chi$ varies over even (non-principal) characters or over odd characters. For example, for each fixed $q \in \N$ Munsch and Shparlinski~\cite{munshpar} proved conjecturally sharp lower bounds for the $2q$-th moments, namely
$$ \frac{1}{r-1} \sum_{\substack{\chi \; \text{mod} \; r, \\ \chi \; \text{even}, \\ \chi \neq \chi_0}} |\theta(1, \chi)|^{2q} \gg_{q} r^{q/2} \log^{(q-1)^2}r , \;\;\;\;\;\;\;\; \frac{1}{r-1} \sum_{\substack{\chi \; \text{mod} \; r, \\ \chi \; \text{odd}}} |\theta(1, \chi)|^{2q} \gg_{q} r^{3q/2} \log^{(q-1)^2}r . $$
Munsch~\cite{munsch} proved almost sharp upper bounds assuming the Generalised Riemann Hypothesis for Dirichlet $L$-functions (losing a factor $\log^{\epsilon}r$ compared with the presumed truth), again when $q \in \N$. One application of such estimates is deducing non-vanishing results for $\theta(1, \chi)$, the subject of a conjecture of Louboutin~\cite{louboutin}. Thus Louboutin~\cite{louboutinodd}, and Louboutin and Munsch~\cite{loubmun}, proved that $\theta(1,\chi) \neq 0$ for $\gg r/\log r$ odd characters and $\gg r/\log r$ even characters modulo prime $r$, by computing and comparing second and fourth moments. Most recently, La Bret\`eche, Munsch and Tenenbaum~\cite{bretechemunschten} introduced weights coming from G\'al sums into such arguments, and proved that $\theta(1,\chi) \neq 0$ for $\gg r/\log^{\delta + o(1)}r$ even characters modulo prime $r$, for a certain explicit constant $\delta \approx 0.086$. See also the work of Bengoechea~\cite{bengo} and of Guo and Peng~\cite{guopeng}, who use Galois theoretic techniques to deduce that $\theta(1,\chi) \neq 0$ for almost all $\chi$, but only for moduli $r$ from certain sparse subsets of primes.

Using our methods, we can prove:

\begin{cor2}\label{cortheta}
Let $r$ be a large prime. Uniformly for any $0 \leq q \leq 1$, we have
$$ \frac{1}{r-1} \sum_{\substack{\chi \; \text{mod} \; r, \\ \chi \; \text{even}}} |\theta(1, \chi)|^{2q} \ll \Biggl(\frac{\sqrt{r}}{1+(1-q)\sqrt{\log\log r}} \Biggr)^q , $$
and
$$ \frac{1}{r-1} \sum_{\substack{\chi \; \text{mod} \; r, \\ \chi \; \text{odd}}} |\theta(1, \chi)|^{2q} \ll \Biggl(\frac{r^{3/2}}{1+(1-q)\sqrt{\log\log r}} \Biggr)^q . $$
\end{cor2}

Since we have $\theta(1, \chi) = \sum_{n=1}^{\infty} \chi(n) e^{-\pi n^{2}/r} \approx \sum_{n \leq \sqrt{r}} \chi(n)$ for even characters, and $\theta(1, \chi) \approx \sum_{n \leq \sqrt{r}} n \chi(n) \approx \sqrt{r} \sum_{n \leq \sqrt{r}} \chi(n)$ for odd characters, one sees that Corollary \ref{cortheta} should follow fairly directly from Theorem \ref{mainthmchar}. Again, see section \ref{secotherproofs} below for the full proof, which is just a partial summation argument.

Similarly as for Theorems \ref{mainthmchar} and \ref{mainthmcont}, it is reasonable to conjecture that one should have {\em matching lower bounds} for the moments in Corollary \ref{cortheta}. If one could prove this for any fixed $0 < q < 1$, then (since the low moments scale proportionally with the exponent $q$, unlike the higher moments where there is quadratic dependence in the exponent of $\log r$) another standard Paley--Zygmund type argument with H\"{o}lder's inequality, comparing upper and lower bounds for two low moments, would immediately yield that $\theta(1, \chi) \neq 0$ for a {\em positive proportion} of $\chi$.

\vspace{12pt}
We are also highly interested in variants of Theorems \ref{mainthmchar} and \ref{mainthmcont}, where the sums are more exotic. In the case of the second moment of such sums, one has exactly the same mean value estimates $\frac{1}{r-1} \sum_{\chi \; \text{mod} \; r} |\sum_{n \leq x} a_n \chi(n)|^2 = \lfloor x \rfloor$ (where $x < r$) and $\frac{1}{T} \int_{0}^{T} |\sum_{n \leq x} a_n n^{it}|^2 dt = \lfloor x \rfloor + O(\frac{x^2}{T})$ for {\em any} complex coefficients $a_n$ with absolute value 1. We cannot seek analogues of Theorems \ref{mainthmchar} and \ref{mainthmcont} at this level of generality, since for a generic sequence of unimodular coefficients $a_n$ the moments will have order $x^q$ (e.g. for independent, uniformly random $a_n$, this is implied by standard moment results like Khintchine's inequalities). But if we tweak the sums in ways that don't disrupt the {\em multiplicative structure} too much, then it turns out that better than squareroot bounds like Theorems \ref{mainthmchar} and \ref{mainthmcont} do endure. Thus we have:

\begin{thm3}\label{mainthmmulttwist}
Let $r$ be a large prime. Then uniformly for any $1 \leq x \leq r$, any $0 \leq q \leq 1$, and any multiplicative function $h(n)$ that has absolute value 1 on primes and absolute value at most 1 on prime powers, we have
$$ \frac{1}{r-1} \sum_{\chi \; \text{mod} \; r} |\sum_{n \leq x} h(n) \chi(n)|^{2q} \ll \Biggl(\frac{x}{1 + (1-q)\sqrt{\log\log(10L)}} \Biggr)^q , $$
where $L = L_r := \min\{x,r/x\}$.
\end{thm3}

The size conditions on $h(n)$ could be adjusted a little here, but the current formulation already permits important examples such as the M\"{o}bius function $\mu(n)$. The key point is that the multiplicative structure of $h(n)$ allows one to adapt the proof of Theorem \ref{mainthmchar}, to show that the moments of $\sum_{n \leq x} h(n) \chi(n)$ behave like the moments of twisted random multiplicative sums $\sum_{n \leq x} h(n) f(n)$. And since, on the primes, $h(p)f(p)$ are again independent and uniform on the complex unit circle (here we use the unimodularity of $h(p)$), these are essentially the same as the untwisted moments. See section \ref{secotherproofs} for further details. One could also prove the analogous result for $\sum_{n \leq x} h(n) n^{it}$.

Unlike the sums $\sum_{n \leq x} \chi(n)$ and $\sum_{n \leq x} n^{it}$, there is no symmetry or ``Fourier flip'' in the presence of a general multiplicative twist $h(n)$. So whilst the appearance of $L_r, L_T$ in Theorems \ref{mainthmchar} and \ref{mainthmcont} was very natural, and one would conjecture those bounds to be sharp, for many $h(n)$ it seems likely that Theorem \ref{mainthmmulttwist} should hold in a stronger form with $\log\log(10L)$ replaced by $\log\log(10x)$. Similarly, the restriction to $x \leq r$ or to $x \leq T$ no longer seems natural for $\sum_{n \leq x} h(n) \chi(n)$ and $\sum_{n \leq x} h(n) n^{it}$. Indeed, if we take $h(n)$ to be a Steinhaus random multiplicative function, then \eqref{precisehelson} implies that for {\em any} large $x$ and large prime $r$ we have
$$ \E \frac{1}{r-1} \sum_{\chi \; \text{mod} \; r} |\sum_{n \leq x} h(n) \chi(n)|^{2q} =  \frac{1}{r-1} \sum_{\chi \; \text{mod} \; r} \E |\sum_{n \leq x} h(n) \chi(n)|^{2q} \asymp \Biggl(\frac{x}{1 + (1-q)\sqrt{\log\log x}} \Biggr)^q , $$
since $h(n) \chi(n)$ will also be a Steinhaus random multiplicative function\footnote{Strictly speaking, $h(n) \chi(n)$ will be a Steinhaus random multiplicative function restricted to numbers not divisible by $r$, so we will have $\E |\sum_{n \leq x} h(n) \chi(n)|^{2q} = \E |\sum_{n \leq x} h(n) - h(r) \sum_{n \leq x/r} h(n)|^{2q}$. But it is easy to check this also obeys the bounds \eqref{precisehelson}, for example by very slightly modifying the proofs of Harper~\cite{harperrmflowmoments} (only a single Euler product factor corresponding to the prime $r$ must change), or simply using Minkowski's inequality (when $1/2 \leq q \leq 1$) and H\"{o}lder's inequality.} for any given Dirichlet character $\chi$. Consequently, for most realisations of $h(n)$ we must have the stronger bound $\frac{1}{r-1} \sum_{\chi \; \text{mod} \; r} |\sum_{n \leq x} h(n) \chi(n)|^{2q} \ll (\frac{x}{1 + (1-q)\sqrt{\log\log x}} )^q$.

In particular, we can return to the case of the M\"{o}bius function. One tends to think of $\mu(n)$ as being ``random looking'' in various senses, and in a recent paper Gorodetsky~\cite{gorodetsky} explored this in the context of character sums. Based on function field considerations, he ultimately conjectured that for all natural number exponents $q < \log r$, the moments $\frac{1}{r-1} \sum_{\chi \; \text{mod} \; r} |\sum_{n \leq x} \mu(n) \chi(n)|^{2q}$ should be asymptotic to the corresponding random multiplicative moments as $x$ and $r$ become large. Likewise, it now seems reasonable to conjecture that for all $0 \leq q \leq 1$ and any fixed $A > 0$, we should have
$$ \frac{1}{r-1} \sum_{\chi \; \text{mod} \; r} |\sum_{n \leq x} \mu(n) \chi(n)|^{2q} \ll (\frac{x}{1 + (1-q)\sqrt{\log\log x}} )^q \;\;\; \forall \; x \leq r^A , $$
and
\begin{equation}\label{mobiusconj}
\frac{1}{2T} \int_{-T}^{T} |\sum_{n \leq x} \mu(n) n^{it}|^{2q} dt \ll (\frac{x}{1 + (1-q)\sqrt{\log\log x}} )^q \;\;\; \forall \; x \leq T^A .
\end{equation}
This conjecture seems worthy of further investigation, since if true it would have significant arithmetic consequences. Standard arguments with Perron's formula imply that
$$ \left| \sum_{x < n \leq x+y} \mu(n) \right| \approx \left| \frac{1}{2\pi i} \int_{-i(x/y)}^{i(x/y)} (\sum_{n \leq 2x} \frac{\mu(n)}{n^s}) \frac{x^{s} ((1+y/x)^s - 1)}{s} ds \right| \lesssim \frac{y}{x} \int_{-x/y}^{x/y} |\sum_{n \leq 2x} \frac{\mu(n)}{n^{it}}| dt , $$
and so \eqref{mobiusconj} (if true) would deliver a bound $|\sum_{x < n \leq x+y} \mu(n) | \lesssim \frac{\sqrt{x}}{(\log\log x)^{1/4}}$ provided $y \leq x^{1-\epsilon}$, say. Thus we could deduce there is cancellation in sums of the M\"{o}bius function in {\em all} short intervals of length $y \gg \frac{\sqrt{x}}{(\log\log x)^{1/4}}$. It is a major open problem to go below, or even to reach, the squareroot interval barrier in problems of this kind. See e.g. the recent beautiful work of Matom\"{a}ki and Radziwi\l\l~\cite{matradz}, which (among many other results) establishes the existence of positive and negative M\"{o}bius values (or, strictly speaking, values of the closely related Liouville function) in intervals of length $C\sqrt{x}$, for a large constant $C$.

\subsection{Ideas from the proofs}\label{subsecproofideas}
Sections \ref{secprep} and \ref{secproof1} will be taken up with the proof of Theorem \ref{mainthmchar}. Here we try to explain and motivate the main steps in the argument.

Throughout the paper, $r$ will be a large prime modulus as in Theorem \ref{mainthmchar}, and we shall write $\Echar$ to denote averaging over all Dirichlet characters mod $r$. Thus if $W(\chi)$ is any function, then $\Echar W := \frac{1}{r-1} \sum_{\chi \; \text{mod} \; r} W(\chi)$. This notation is intended both to simplify the writing, and to emphasise the connection between the arguments here and those in the random multiplicative function case. 

\vspace{12pt}
The introduction to the author's paper~\cite{harperrmflowmoments} contains a detailed outline of the proof of the estimate \eqref{precisehelson} for $\E|\sum_{n \leq x} f(n)|^{2q}$. The first phase of that proof involves showing that
$$ \E|\sum_{n \leq x} f(n)|^{2q} \approx x^{q} \E\Biggl(\frac{1}{\log x} \int_{-1/2}^{1/2} |F(1/2+it)|^2 dt \Biggr)^{q} , $$
where $F(s) := \prod_{p \leq x} (1 - \frac{f(p)}{p^{s}})^{-1}$ is a suitable random Euler product. The second phase is to show that $\E (\frac{1}{\log x} \int_{-1/2}^{1/2} |F(1/2+it)|^2 dt )^{q} \approx (\frac{1}{1 + (1-q)\sqrt{\log\log x}} )^q$, using ideas from the probabilistic theory of {\em critical multiplicative chaos}. To prove Theorem \ref{mainthmchar}, we shall rework the first phase of the random multiplicative function proof to show, roughly speaking, that
\begin{equation}\label{charbyranddisplay}
\Echar |\sum_{n \leq x} \chi(n)|^{2q} \ll x^{q} \E\Biggl(\frac{1}{\log P} \int_{-1/2}^{1/2} |F_{P}^{\text{rand}}(1/2+it)|^2 dt \Biggr)^{q} ,
\end{equation}
where $F_{P}^{\text{rand}}(s) := \prod_{p \leq P} (1 - \frac{f(p)}{p^s})^{-1}$ for a suitable parameter $P = P(r,x)$. Notice that on the left hand side we have our character average, but we will show this can be bounded by the expectation of the genuinely random quantity on the right. Thus we will {\em not} need to rework the second phase of the random multiplicative function proof, since we will be able to apply this directly once we have established a bound like \eqref{charbyranddisplay}.

\vspace{12pt}
In \cite{harperrmflowmoments}, to make the connection between $\E|\sum_{n \leq x} f(n)|^{2q}$ and a random Euler product one conditions on (i.e. freezes) the values $f(p)$ for all ``small'' primes $p$, before applying H\"{o}lder's inequality to compare with the conditional second moment. We shall review that argument now, but for simplicity only for $\sum_{n \leq x, P(n) > \sqrt{x}} f(n)$, the subsum over numbers whose largest prime factor $P(n)$ is greater than $\sqrt{x}$.

By multiplicativity we can write $\sum_{n \leq x, P(n) > \sqrt{x}} f(n) = \sum_{\sqrt{x} < p \leq x} f(p) \sum_{m \leq x/p} f(m)$, and note that in the inner sums we have $m \leq x/p < \sqrt{x}$, so those values are entirely determined by the $(f(p))_{p \leq \sqrt{x}}$. Then if we let $\tilde{\E}$ denote expectation conditional on the values $(f(p))_{p \leq \sqrt{x}}$ (i.e. expectation with those values treated as fixed and the $(f(p))_{p > \sqrt{x}}$ remaining random, so the conditional expectation of any quantity is a function of the values $(f(p))_{p \leq \sqrt{x}}$), we get
$$ \tilde{\E} \Biggl| \sum_{n \leq x, P(n) > \sqrt{x}} f(n) \Biggr|^{2} = \sum_{\sqrt{x} < p,q \leq x} \tilde{\E}\Biggl( f(p) (\sum_{m \leq x/p} f(m)) \overline{f(q)} (\sum_{m \leq x/q} \overline{f(m)}) \Biggr) = \sum_{\sqrt{x} < p \leq x} \Biggl| \sum_{m \leq x/p} f(m) \Biggr|^2 , $$
because the $(f(p))_{p > \sqrt{x}}$ {\em remain independent with mean zero}. But by the Tower Property of conditional expectation (which in this case is simply Fubini's theorem, breaking up the multiple ``integration'' $\E$ into separate integrations corresponding to the $(f(p))_{p \leq \sqrt{x}}$ and the $(f(p))_{p > \sqrt{x}}$), we can write $\E|\sum_{n \leq x, P(n) > \sqrt{x}} f(n) |^{2q} = \E \tilde{\E} | \sum_{n \leq x, P(n) > \sqrt{x}} f(n) |^{2q}$. So applying H\"{o}lder's inequality to the {\em conditional expectation} $\tilde{\E}$, we get
$$ \E\Biggl|\sum_{n \leq x, P(n) > \sqrt{x}} f(n) \Biggr|^{2q} \leq \E \Biggl( \tilde{\E} \Biggl| \sum_{n \leq x, P(n) > \sqrt{x}} f(n) \Biggr|^{2} \Biggr)^{q} = \E \Biggl( \sum_{\sqrt{x} < p \leq x} \Biggl| \sum_{m \leq x/p} f(m) \Biggr|^2 \Biggr)^{q} . $$
Having reached this point, one can perform some smoothing and use a suitable form of Parseval's identity to relate the right hand side to an Euler product average like $\E(\frac{1}{\log x} \int_{-1/2}^{1/2} |F(1/2+it)|^2 dt )^{q}$. The factor $\frac{1}{\log x}$ reflects the density of primes on the interval $(\sqrt{x},x]$.

To prove Theorem \ref{mainthmchar}, we must find a version of the above procedure that can succeed for character averages, and ultimately deliver \eqref{charbyranddisplay}. This is challenging for two related reasons: when averaging over characters $\chi$ mod $r$, the values $\chi(p)$ for different primes $p$ are not really independent of one another; and since there are only $r-1$ characters, we cannot hope for the behaviour of the $\chi(n)$ to really resemble the random model $f(n)$ on events that occur with probability much smaller than $1/(r-1)$ on the random side. In particular, if we try to condition (freeze) the values $\chi(p)$ for too many primes $p$, there will generally be zero or precisely one character with the values $\chi(p)$ as specified, and we cannot hope for conditioning on such configurations to match the calculations in the random multiplicative case.

\vspace{12pt}
To match character averages with expectations involving random multiplicative functions, we must confine ourselves to working with {\em polynomial expressions of suitable degree} in the various $\chi(p)$. Indeed, if $p_1, ..., p_k, p_{k+1}, ..., p_l$ are any (not necessarily distinct) primes such that $\prod_{j=1}^{k} p_j, \prod_{j=k+1}^{l} p_j < r$, we have a genuine equality
\begin{eqnarray}\label{charrmfexp}
\Echar \prod_{j=1}^{k} \chi(p_j) \prod_{j=k+1}^{l} \overline{\chi(p_j)} & = & \Echar \chi(\prod_{j=1}^{k} p_j) \overline{\chi(\prod_{j=k+1}^{l} p_j)} = \textbf{1}_{\prod_{j=1}^{k} p_j \equiv \prod_{j=k+1}^{l} p_j \; \text{mod} \; r} \nonumber \\
& = & \textbf{1}_{\prod_{j=1}^{k} p_j = \prod_{j=k+1}^{l} p_j} = \E \prod_{j=1}^{k} f(p_j) \prod_{j=k+1}^{l} \overline{f(p_j)} .
\end{eqnarray}
Here we used orthogonality of Dirichlet characters; the size restrictions on $\prod_{j=1}^{k} p_j$ and $\prod_{j=k+1}^{l} p_j$ (to crucially replace a congruence mod $r$ by an equality); and finally the orthogonality property $\E f(n) \overline{f(m)} = \textbf{1}_{n=m}$ of Steinhaus random multiplicative functions. 

To exploit \eqref{charrmfexp}, the key observation is that we did not really need to condition on the exact values of all $(f(p))_{p \leq \sqrt{x}}$ to carry our overall argument through on the random multiplicative side. Since we ultimately bound $\sum_{\sqrt{x} < p \leq x} | \sum_{m \leq x/p} f(m) |^2$ by something like $\frac{1}{\log x} \int_{-1/2}^{1/2} |F(1/2+it)|^2 dt$, the proof could proceed very similarly if we conditioned on any information about the $f(p)$ that is sufficient to roughly fix the value of $\int_{-1/2}^{1/2} |F(1/2+it)|^2 dt$. For each $t$, the value of $|F(1/2+it)|$ is the exponential of a prime number sum like $\Re \sum_{p \leq \sqrt{x}} \frac{f(p)}{p^{1/2+it}}$. Furthermore, such prime number sums do not usually vary much when $t$ varies by less than $1/\log x$ (say), because $p^{it} = e^{it\log p}$ does not vary much when $t$ varies by much less than $1/\log p$. So if we replace $\tilde{\E}$ in the above description by a {\em coarser conditioning}, only on the values of $\Re \sum_{p \leq \sqrt{x}} \frac{f(p)}{p^{1/2+it}}$ at some net of $t$ values with spacing roughly $1/\log x$, (and in fact we only need to know $\Re \sum_{p \leq \sqrt{x}} \frac{f(p)}{p^{1/2+it}}$ to a precision of $O(1)$), this should fix enough information on the inside of the $q$-th power that we still end up with an overall bound of roughly the shape $\E(\frac{1}{\log x} \int_{-1/2}^{1/2} |F(1/2+it)|^2 dt )^{q}$.

Now translating to character sums, the analogue of breaking up $\E|\sum_{n \leq x, P(n) > \sqrt{x}} f(n) |^{2q}$ by writing $\E|\sum_{n \leq x, P(n) > \sqrt{x}} f(n) |^{2q} = \E \tilde{\E} | \sum_{n \leq x, P(n) > \sqrt{x}} f(n) |^{2q}$ will be (roughly) to write
$$ \Echar |\sum_{n \leq x} \chi(n)|^{2q} = \sum_{\textbf{j}} \Echar G_{\textbf{j}}(\chi) |\sum_{n \leq x} \chi(n) |^{2q} , $$
where $G_{\textbf{j}}(\chi)$ are functions satisfying $\sum_{\textbf{j}} G_{\textbf{j}}(\chi) \equiv 1$ that approximately pick out all characters $\chi$ for which sums like $\Re \sum_{p \leq \sqrt{x}} \frac{\chi(p)}{p^{1/2+it}}$ have a given collection of values. We can apply H\"{o}lder's inequality to each inner average $\Echar G_{\textbf{j}}(\chi) |\sum_{n \leq x} \chi(n) |^{2q}$ here, as we did with $\tilde{\E} | \sum_{n \leq x, P(n) > \sqrt{x}} f(n) |^{2q}$ previously, and end up needing to work with quantities like $\Echar G_{\textbf{j}}(\chi) |\sum_{n \leq x} \chi(n) |^{2}$. Furthermore, if the functions $G_{\textbf{j}}$ are sufficiently nice we could hope to approximate $G_{\textbf{j}}(\chi)$ by a polynomial in the $\chi(p)$, at which point we might be able to invoke \eqref{charrmfexp} and replace $\Echar G_{\textbf{j}}(\chi) |\sum_{n \leq x} \chi(n) |^{2}$ by $\E G_{\textbf{j}}(f) |\sum_{n \leq x} f(n) |^{2}$. This is a fairly good high level description of how the proof of Theorem \ref{mainthmchar} proceeds, and so a key tool will be a collection of nice functions forming a smooth partition of unity, from which we can ultimately form our multipliers $G_{\textbf{j}}$ by taking suitable products (corresponding to all the different $t$ values in our net). See Approximation Result \ref{apres}, in section \ref{subsecpartition}, for the technical statement we use.

\vspace{12pt}
But there is an important issue that must be addressed. Notice that the conductor $r$ doesn't appear in the preceding paragraph, and one doesn't yet see how the quantity $\log\log(10L_r)$ will arise in Theorem \ref{mainthmchar} in place of $\log\log x$. The point is that in \eqref{charrmfexp} we must ensure, when everything is expanded out, that our products of primes are smaller than $r$. See Propositions \ref{condexpprop} and \ref{boxprobprop}, in section \ref{subsecmeans}, for the precise statements that we will use to compare $\Echar G_{\textbf{j}}(\chi) |\sum_{n \leq x} \chi(n) |^{2}$ with $\E G_{\textbf{j}}(f) |\sum_{n \leq x} f(n) |^{2}$, where it is crucial that the prime sums $\Re \sum_{p \leq P} \frac{\chi(p)}{p^{1/2+it}}$ involved run up to $P$ for some $P$ that is suitably small compared with $r/x$. Indeed, for prime number sums of length $P$ we should work with a net of $t$ values with spacing roughly $1/\log P$, so with roughly $\log P$ different sums. And for each of those, when we apply Taylor's theorem to expand its contribution to $G_{\textbf{j}}(\chi)$ we must expand as far as degree $\log^{O(1)}P$ to achieve a suitable level of precision. Since the square $|\sum_{n \leq x} \chi(n) |^{2}$ also contributes terms $\chi(n), \overline{\chi(m)}$ with $n,m$ up to $x$ in size, we must have $x P^{\log^{O(1)}P} < r$ to match up character sum and random multiplicative function expressions.

In particular, we cannot work with $P=\sqrt{x}$ here unless $x$ is rather small compared with $r$. The largest permissible choice of $P$ will rather be something like $e^{\log^{c}(r/x)}$, for a suitable constant $c$. Since we also want $P < x$ (it would not make sense to include sums over primes $> x$, which are not involved in $\sum_{n \leq x} \chi(n)$, in the conditioning), one reasonable choice turns out to be $P \approx \exp\{\log^{1/6}L_r\}$. Fortunately, for an upper bound (but not for a lower bound) we are free to select the range of primes involved in our ``conditioning'' as we wish. (This observation would actually allow some simplification of the upper bound arguments in the purely random setting~\cite{harperrmflowmoments} as well. Rather than breaking up $\sum_{n \leq x} f(n)$ into various subsums according to the size of the largest prime factor $P(n)$, as there, and then conditioning on all somewhat smaller primes, we can simply condition the full sum on all primes up to one suitably chosen point.) As $P$ becomes smaller, the possible saving $(\frac{1}{1 + (1-q)\sqrt{\log\log P}} )^q$ that we can ultimately achieve using multiplicative chaos results will diminish, but since $\log\log P$ varies so slowly we can vary $P$ quite a lot without visibly changing the final bounds. In particular, note that although our choice $P \approx \exp\{\log^{1/6}L_r\}$ is (necessarily) significantly smaller than $L_r$, we still have $\log\log P \asymp \log\log L_r$ and so we get the desired factor $\frac{1}{1 + (1-q)\sqrt{\log\log(10L)}}$ in Theorem \ref{mainthmchar}.

In sections \ref{thm1small}--\ref{subsecmaingorandom}, we implement an argument along the above lines and succeed in replacing all terms of the form $\Echar G_{\textbf{j}}(\chi) |\sum_{n \leq x} \chi(n) |^{2}$ by $\E G_{\textbf{j}}(f) |\sum_{n \leq x} f(n) |^{2}$, thus passing to the purely random setting. However, because this is done in the setup of coarser conditioning, more work remains to confirm that the resulting random multiplicative function expression can really still be bounded by something like $\E(\frac{1}{\log P} \int_{-1/2}^{1/2} |F_{P}^{\text{rand}}(1/2+it)|^2 dt )^{q}$. Initially one can apply a similar kind of smoothing and Parseval procedure as was originally done in \cite{harperrmflowmoments}, see section \ref{subsecparseval} (where we recall a version of Parseval's identity for Dirichlet series) and section \ref{subsecpasseuler} below. This brings us to random Euler products, but with a (weighted) sum over $\textbf{j}$ on the outside of the $q$-th power in place of a ``genuine'' expectation $\E$, and with some additional averaging (weighted by $G_{\textbf{j}}(f)$) on the inside of the $q$-th power. Since each term $G_{\textbf{j}}(f)$ only contains information about $\Re \sum_{p \leq P} \frac{f(p)}{p^{1/2+it}}$ at a net of points $t$, as opposed to all $-1/2 \leq t \leq 1/2$ or all $t \in \R$, we want to restrict to working with Euler products at those $t$. However, provided the net of $t$ are sufficiently close together (slightly closer than $1/\log P$), this can be achieved with simple mean square arguments.

Finally, we must check that the weighting by $G_{\textbf{j}}(f)$ fixes enough information about all the sums $\Re \sum_{p \leq P} \frac{f(p)}{p^{1/2+it}}$ that, even with this extra averaging on the inside of the $q$-th power, we still end up with something roughly like $\E(\frac{1}{\log P} \int_{-1/2}^{1/2} |F_{P}^{\text{rand}}(1/2+it)|^2 dt )^{q}$ (or actually a version of this with the integral replaced by a sum over our discrete points $t$). Provided the parameters were selected properly in terms of $P$ when constructing the smooth functions underlying $G_{\textbf{j}}(\cdot)$, (which is ultimately why one must Taylor expand as far as degree $\log^{O(1)}P$ in the above discussion), it turns out that averaging against $G_{\textbf{j}}(f)$ does essentially restrict all the sums $\Re \sum_{p \leq P} \frac{f(p)}{p^{1/2+it}}$ to their intended boxes depending on $\textbf{j}$. Thus the weighted outer sum over $\textbf{j}$ does perform essentially the same role as a genuine expectation, and we do arrive at (a discretised version of) $\E(\frac{1}{\log P} \int_{-1/2}^{1/2} |F_{P}^{\text{rand}}(1/2+it)|^2 dt )^{q}$. This argument is completed in sections \ref{subsecrefinecond}--\ref{subsecconclusion}, using properties of the functions from Approximation Result \ref{apres} and using multiplicative chaos results presented in section \ref{subsecrandeuler}.

\section{Preparations}\label{secprep}

\subsection{A smooth partition of unity}\label{subsecpartition}

As explained in the introduction, one important tool for us will be a collection of fairly well behaved functions that we can use to approximately detect the values of various sums involving $\chi(p)$, and therefore simulate a conditioning process in our main proofs. These sorts of constructions are fairly standard in modern analysis and number theory, and it will not be too difficult to prove the following.

\begin{approxres1}\label{apres}
Let $N \in \N$ be large, and $\delta > 0$ be small. There exist functions $g : \R \rightarrow \R$ (depending on $\delta$) and $g_{N+1} : \R \rightarrow \R$ (depending on $\delta$ and $N$) such that, if we define $g_{j}(x) = g(x - j)$ for all integers $|j| \leq N$, we have the following properties:
\begin{enumerate}
\item $\sum_{|j| \leq N} g_{j}(x) + g_{N+1}(x) = 1$ for all $x \in \R$;

\item $g(x) \geq 0$ for all $x \in \R$, and $g(x) \leq \delta$ whenever $|x| > 1$;

\item $g_{N+1}(x) \geq 0$ for all $x \in \R$, and $g_{N+1}(x) \leq \delta$ whenever $|x| \leq N$;

\item for all $l \in \N$ and all $x \in \R$, we have the derivative estimate $|\frac{d^{l}}{dx^{l}} g(x)| \leq \frac{1}{\pi (l+1)} (\frac{2\pi}{\delta})^{l+1}$.

\end{enumerate}
\end{approxres1}

\begin{proof}[Proof of Approximation Result \ref{apres}]
Our construction will be a minor variant of the proof of Approximation Result 1 of Harper~\cite{harperpartition} (with $R=0$ there). As such, we content ourselves with outlining the main steps. 

Let $b(x)$ be a Beurling--Selberg function majorising the indicator function $\textbf{1}_{|x| \leq 1/2}$, with Fourier transform supported on $[-1/\delta,1/\delta]$. See e.g. Vaaler's paper~\cite{vaaler} for background on such majorants. Thus we have $b(x) \geq \textbf{1}_{|x| \leq 1/2}$ for all $x \in \R$; and $\int_{-\infty}^{\infty} b(x) dx = 1 + \delta$; and $b(x) = \int_{-1/\delta}^{1/\delta} \hat{b}(t) e^{2\pi i xt} dt$ for all $x \in \R$, where $|\hat{b}(t)| = |\int b(x) e^{-2\pi i x t} dx| \leq 1 + \delta$.

We define $g(x)$ as a convolution of $b$, namely
$$ g(x) = \int_{-\infty}^{\infty} \textbf{1}_{|u| \leq 1/2} \frac{b(x-u)}{1+\delta} du = \int_{-\infty}^{\infty} \textbf{1}_{|x-u| \leq 1/2} \frac{b(u)}{1+\delta} du . $$
Then it is clear that $g(x)$ is non-negative, since $b(x)$ is non-negative. The other claims about $g(x)$ in (ii) and (iv) follow identically as in Harper~\cite{harperpartition}.

Now by definition of $g$ and $g_j$, the sum $\sum_{|j| \leq N} g_{j}(x)$ is
$$ = \int_{-\infty}^{\infty} \sum_{|j| \leq N} \textbf{1}_{|x - j -u| \leq 1/2} \frac{b(u)}{1 + \delta} du = \int_{-\infty}^{\infty} \textbf{1}_{|x-u| \leq N+1/2} \frac{b(u)}{1 + \delta} du $$
for all real $x$. Thus we always have $\sum_{|j| \leq N} g_{j}(x) \leq \int_{-\infty}^{\infty} \frac{b(u)}{1 + \delta} du = 1$, and furthermore $1 - \sum_{|j| \leq N} g_{j}(x)$ is equal to
$$ \int_{-\infty}^{\infty} \textbf{1}_{|x-u| > N+1/2} \frac{b(u)}{1 + \delta} du = \int_{-\infty}^{\infty} \textbf{1}_{|x-u| > N+1/2} \frac{b(u) - \textbf{1}_{|u| \leq 1/2}}{1 + \delta} du + \int_{-\infty}^{\infty} \textbf{1}_{|x-u| > N+1/2} \frac{\textbf{1}_{|u| \leq 1/2}}{1 + \delta} du . $$
When $|x| \leq N$ the second integral here vanishes, and so
$$ 1 - \sum_{|j| \leq N} g_{j}(x) \leq \int_{-\infty}^{\infty} \frac{|b(u) - \textbf{1}_{|u| \leq 1/2}|}{1+\delta} du = \int_{-\infty}^{\infty} \frac{b(u) - \textbf{1}_{|u| \leq 1/2}}{1+\delta} du = \frac{\delta}{1+\delta} \leq \delta . $$
So the first and third statements in Approximation Result \ref{apres} follow if we simply set $g_{N+1}(x):= 1 - \sum_{|j| \leq N} g_{j}(x)$ for all $x \in \R$.
\end{proof}

\subsection{Mean value estimates}\label{subsecmeans}

Having introduced approximating functions as in Approximation Result \ref{apres}, our arguments will require us to evaluate various character averages involving these functions. A basic tool for this will be the following bound, which is a fairly standard even moment estimate for character sums of appropriate lengths.

\begin{lem1}[Even moment estimate]\label{evenmomentlem}
Let $x \geq 1$, and let $(c(n))_{n \leq x}$ be any complex numbers. Let $\mathcal{P}$ be any finite set of primes, let $\mathcal{Q}$ be any (non-empty) set consisting of some elements of $\mathcal{P}$ and squares of elements of $\mathcal{P}$, and write $U := \max\{q \in \mathcal{Q}\}$ . Finally, let $Q(\chi) := \sum_{q \in \mathcal{Q}} \frac{a(q) \chi(q)}{\sqrt{q}}$, where the $a(q)$ are any complex numbers.

Then for any natural number $k$ such that $x U^k < r$, we have
$$ \Echar \Biggl|\sum_{n \leq x} c(n) \chi(n) \Biggr|^{2} |Q(\chi)|^{2k} \ll \Biggl(\sum_{n \leq x} \tilde{d}(n) |c(n)|^2 \Biggr) \cdot (k !) \Biggl( 2 \sum_{q \in \mathcal{Q}} \frac{v_q |a(q)|^2}{q} \Biggr)^{k} , $$
where $\tilde{d}(n) := \sum_{d|n} \textbf{1}_{p|d \Rightarrow p \in \mathcal{P}}$, and $v_q$ is 1 if $q$ is a prime and 6 if $q$ is the square of a prime.
\end{lem1}

\begin{proof}[Proof of Lemma \ref{evenmomentlem}]
This is a character sum version of Lemma 2 of Harper~\cite{harperpartition}, which dealt with $t$-averages of $\sum_{n \leq x} c(n) n^{-it}$ and $\sum_{q \in \mathcal{Q}} \frac{a(q)}{q^{1/2 + it}}$. It may be proved in the same way as that result (in fact slightly more easily, since for Dirichlet characters one has perfect rather than approximate orthogonality).
\end{proof}

Using Taylor expansion and Lemma \ref{evenmomentlem} we can deduce the following crucial Proposition, which guarantees that (provided we keep sufficient control on the sizes of the various parameters) character averages involving the functions $g_j$ behave in the same way as the corresponding averages involving random multiplicative functions.

\begin{prop1}[Characters behave like random model]\label{condexpprop}
Let the functions $g_j$, with associated parameters $N, \delta$, be as in Approximation Result \ref{apres}. Suppose that $x \geq 1$, and let $(c(n))_{n \leq x}$ be any complex numbers having absolute values $\leq 1$. Furthermore, let $P$ be large, and let $Y \in \N$ be such that $x P^{400(Y/\delta)^2 \log(N\log P)} < r$. Let $f(n)$ denote a Steinhaus random multiplicative function.

Then for any indices $-N \leq j(1), j(2), ..., j(Y) \leq N+1$, and for any sequences $(a_{1}(p))_{p \leq P}, (a_{1}(p^2))_{p \leq P}, ..., (a_{Y}(p))_{p \leq P}, (a_{Y}(p^2))_{p \leq P}$ of complex numbers having absolute values $\leq 1$, we have
\begin{eqnarray}
&& \Echar \prod_{i=1}^{Y} g_{j(i)}(\Re(\sum_{p \leq P} \frac{a_{i}(p) \chi(p)}{\sqrt{p}} + \frac{a_{i}(p^2) \chi(p^2)}{p})) \Biggl|\sum_{n \leq x} c(n) \chi(n) \Biggr|^2 \nonumber \\
& = & \E \prod_{i=1}^{Y} g_{j(i)}(\Re(\sum_{p \leq P} \frac{a_{i}(p) f(p)}{\sqrt{p}} + \frac{a_{i}(p^2) f(p^2)}{p})) \Biggl|\sum_{n \leq x} c(n) f(n) \Biggr|^2 + O\left(\frac{x}{(N \log P)^{Y(1/\delta)^2}} \right) . \nonumber
\end{eqnarray}
\end{prop1}

\begin{proof}[Proof of Proposition \ref{condexpprop}]
Using property (iv) from Approximation Result \ref{apres}, for any threshold $2S \in 2\N$ we can write $g_{j}(x) = \tilde{g}_{j}(x) + r_{j}(x)$, where $\tilde{g}_{j}(\cdot)$ is a polynomial of degree $2S-1$ (namely the degree $2S-1$ Taylor polynomial of $g_{j}(x)$ about zero), and where $|r_{j}(x)| \leq \frac{|x|^{2S}}{(2S)!} \sup_{|y| \leq |x|} |\frac{d^{2S}}{dy^{2S}} g_{j}(y)| \ll \frac{N |2\pi x/\delta|^{2S}}{\delta S (2S)!}$. (Note that the factor $N$ here is to account for the case where $g_{j}(y) = g_{N+1}(y) = 1 - \sum_{|i| \leq N} g_{i}(y)$.) First we examine the contribution from the ``main terms'' $\tilde{g}_{j}(\cdot)$ to the left hand side in the Proposition. Provided that $x P^{4SY} < r$, we can expand all the polynomials and the square out and find that
\begin{eqnarray}
&& \Echar \prod_{i=1}^{Y} \tilde{g}_{j(i)}(\Re(\sum_{p \leq P} \frac{a_{i}(p) \chi(p)}{\sqrt{p}} + \frac{a_{i}(p^2) \chi(p^2)}{p})) \Biggl|\sum_{n \leq x} c(n) \chi(n) \Biggr|^2 \nonumber \\
& = & \E \prod_{i=1}^{Y} \tilde{g}_{j(i)}(\Re( \sum_{p \leq P} \frac{a_{i}(p) f(p)}{\sqrt{p}} + \frac{a_{i}(p^2) f(p^2)}{p})) \Biggl|\sum_{n \leq x} c(n) f(n) \Biggr|^2 , \nonumber
\end{eqnarray}
since if $1 \leq U,V < r$ then we have the equality $\Echar \chi(U) \overline{\chi(V)} = \textbf{1}_{U=V} = \E f(U) \overline{f(V)}$. (Note that $\E f(U) \overline{f(V)} = \textbf{1}_{U=V}$ for {\em all} natural numbers $U,V$ on the random multiplicative function side, but on the character side one needs the restriction that $1 \leq U,V < r$ to boost a congruence mod $r$ to an equality.)

Next, dividing up according to the smallest index $i$ at which we get a remainder, we see the contribution from all of the remainders $r_{j(i)}(\cdot)$ to the left hand side in Proposition \ref{condexpprop} is
\begin{eqnarray}
& \ll & \Echar \sum_{i=1}^{Y} \frac{N |2\pi/\delta|^{2S}}{\delta S (2S)!} |\sum_{p \leq P} \frac{a_{i}(p) \chi(p)}{\sqrt{p}} + \frac{a_{i}(p^2) \chi(p^2)}{p}|^{2S} \cdot \nonumber \\
&& \cdot \prod_{l=1}^{i-1} \Biggl(1 + O(\frac{N |2\pi/\delta|^{2S}}{\delta S (2S)!} |\sum_{p \leq P} \frac{a_{l}(p) \chi(p)}{\sqrt{p}} + \frac{a_{l}(p^2) \chi(p^2)}{p}|^{2S})\Biggr) \Biggl|\sum_{n \leq x} c(n) \chi(n) \Biggr|^2 . \nonumber
\end{eqnarray}
Here we used the fact that $|\tilde{g}_{j(l)}(x)| \leq |g_{j(l)}(x)| + |r_{j(l)}(x)| \leq 1 + O(\frac{N |2\pi x/\delta|^{2S}}{\delta S (2S)!})$. Using Lemma \ref{evenmomentlem} and the condition $x P^{4SY} < r$ again, along with the bounds $(jS)! \ll \sqrt{jS} (jS/e)^{jS}$ and $(2S)! \geq (2S/e)^{2S}$ to handle the factorials produced by that lemma, we find this is all
\begin{eqnarray}
& \ll & \Biggl(\sum_{n \leq x} \tilde{d}(n) |c(n)|^2 \Biggr) \cdot \sum_{i=1}^{Y} \frac{N|2\pi/\delta|^{2S}}{\delta S (2S)!} \Biggl( \sqrt{iS} (\frac{iS}{e})^S (2\log\log P + O(1))^{S} \Biggr) \cdot \nonumber \\
&& \cdot \prod_{l=1}^{i-1} \Biggl(1 + O\Biggl(\frac{N|2\pi/\delta|^{2S}}{\delta S (2S)!} (\frac{iS}{e})^S (2\log\log P + O(1))^{S} \Biggr) \Biggr) \nonumber \\
& \ll & x\log P \cdot \sum_{i=1}^{Y} \frac{N}{\delta} \sqrt{\frac{i}{S}} (\frac{e (\pi/\delta)^2 i (2\log\log P + O(1))}{S})^S \cdot \nonumber \\
&& \cdot \Biggl( 1 + O\Biggl(\frac{N}{\delta S} (\frac{e (\pi/\delta)^2 i (2\log\log P + O(1))}{S})^S \Biggr) \Biggr)^{i-1} . \nonumber
\end{eqnarray}
Here we also used our assumptions that $|c(n)| \leq 1$ and $|a_{i}(p)|, |a_{i}(p^2)| \leq 1$, which in particular imply that $2\sum_{p \leq P} (\frac{|a_{i}(p)|^2 }{p} + \frac{6|a_{i}(p^2)|^2}{p^2}) \leq 2\sum_{p \leq P} \frac{1}{p} + O(1) = 2\log\log P + O(1)$. One has the same overall bound for the contribution from the remainders $r_{j(i)}(\cdot)$ to the right hand side in Proposition \ref{condexpprop}, since one has the same bound for $\E |\sum_{n \leq x} c(n) f(n) |^{2} |Q(f)|^{2k}$ as for the character average in Lemma \ref{evenmomentlem} (indeed this quantity is again exactly equal to $\Echar |\sum_{n \leq x} c(n) \chi(n) |^{2} |Q(\chi)|^{2k}$, under the size conditions in the lemma).

Now if we set $S = 100Y \lfloor (1/\delta)^2 \log(N\log P) \rfloor$, then the condition $x P^{4SY} < r$ is satisfied in view of our assumption that $x P^{400(Y/\delta)^2 \log(N\log P)} < r$. And with this choice, the error term produced by all of the remainders is
$$ \ll x\log P \cdot \sum_{i=1}^{Y} \frac{N}{\delta} \sqrt{\frac{i}{S}} 0.6^S \left( 1 + O(\frac{N}{\delta S} 0.6^S ) \right)^{i-1} \ll x\log P \cdot NY \cdot 0.6^S \ll \frac{x}{(N \log P)^{Y(1/\delta)^2}} , $$
as desired.
\end{proof}

Taking $x=1$ and $c(1)=1$ in Proposition \ref{condexpprop}, we obtain the following important special case.

\begin{prop2}\label{boxprobprop}
Let the functions $g_j$, with associated parameters $N, \delta$, be as in Approximation Result \ref{apres}. Suppose $P$ is large, and let $Y \in \N$ be such that $P^{400(Y/\delta)^2 \log(N\log P)} < r$. Let $f(n)$ denote a Steinhaus random multiplicative function.

Then for any indices $-N \leq j(1), j(2), ..., j(Y) \leq N+1$, and for any sequences $(a_{1}(p))_{p \leq P}, (a_{1}(p^2))_{p \leq P}, ..., (a_{Y}(p))_{p \leq P}, (a_{Y}(p^2))_{p \leq P}$ of complex numbers having absolute values $\leq 1$, we have
\begin{eqnarray}
\Echar \prod_{i=1}^{Y} g_{j(i)}(\Re(\sum_{p \leq P} \frac{a_{i}(p) \chi(p)}{\sqrt{p}} + \frac{a_{i}(p^2) \chi(p^2)}{p})) & = & \E \prod_{i=1}^{Y} g_{j(i)}(\Re(\sum_{p \leq P} \frac{a_{i}(p) f(p)}{\sqrt{p}} + \frac{a_{i}(p^2) f(p^2)}{p})) + \nonumber \\
&& + O\left(\frac{1}{(N \log P)^{Y(1/\delta)^2}} \right) . \nonumber
\end{eqnarray}
\end{prop2}

\subsection{Parseval's identity for Dirichlet series}\label{subsecparseval}
As in the proof of the random analogue of Theorem \ref{mainthmchar}, we shall need the following version of Parseval's identity for Dirichlet series.
\begin{harman1}[See (5.26) in sec. 5.1 of Montgomery and Vaughan~\cite{mv}]\label{harmdirichlet}
Let $(a_n)_{n=1}^{\infty}$ be any sequence of complex numbers, and let $A(s) := \sum_{n=1}^{\infty} \frac{a_n}{n^s}$ denote the corresponding Dirichlet series, and $\sigma_c$ denote its abscissa of convergence. Then for any $\sigma > \max\{0,\sigma_c \}$, we have
$$ \int_{0}^{\infty} \frac{|\sum_{n \leq x} a_n |^2}{x^{1 + 2\sigma}} dx = \frac{1}{2\pi} \int_{-\infty}^{\infty} \left|\frac{A(\sigma + it)}{\sigma + it}\right|^2 dt . $$
\end{harman1}

We shall deploy Harmonic Analysis Result \ref{harmdirichlet} at a point in our argument where we have already used Propositions \ref{condexpprop} and \ref{boxprobprop} to move from studying character averages to studying random multiplicative functions $f(n)$. As such, we will be able to take $a_n$ as the values $f(n)$ restricted to $P$-smooth numbers $n$ (for a certain parameter $P$), and with $A(s)$ being the partial Euler product corresponding to $f(n)$ on all $P$-smooth numbers, similarly as in the random case in \cite{harperrmflowmoments}. Note that we could not proceed in this way working with Dirichlet characters mod $r$ directly, because we would not retain sufficient control on terms in the Euler product corresponding to those $n > r$.

\subsection{Random Euler products}\label{subsecrandeuler}
The final key tool we require, and the ultimate source of the better than squareroot cancellation that we look to establish in our theorems, is an upper bound for the small moments of a short integral of a random Euler product. This builds on ideas from the probabilistic theory of critical multiplicative chaos. Recall that $(f(p))_{p \; \text{prime}}$ is a sequence of independent random variables distributed uniformly on the complex unit circle, and for any large quantity $P$ and any complex $s$ with $\Re(s) > 0$, let $F_{P}^{\text{rand}}(s) := \prod_{p \leq P} (1 - \frac{f(p)}{p^s})^{-1}$.

\begin{multchaos1}[See section 4 of Harper~\cite{harperrmflowmoments}]\label{mcres1}
Uniformly for all large $P$ and $2/3 \leq q \leq 1$, we have
$$ \E( \int_{-1/2}^{1/2} |F_{P}^{\text{rand}}(1/2 + it)|^2 dt )^q \ll \Biggl(\frac{\log P}{1 + (1-q)\sqrt{\log\log P}}\Biggr)^{q} . $$
\end{multchaos1}

This is proved in section 4 of \cite{harperrmflowmoments} (see the proof of the Theorem 1 upper bound there), in slightly different notation: the product $F_{0}(1/2+it)$ from \cite{harperrmflowmoments} corresponds to $F_{P}^{\text{rand}}(1/2 + it)$ with $P$ replaced by $x^{1/e}$.

Although the reader is welcome to treat Multiplicative Chaos Result \ref{mcres1} entirely as a black box for our purposes here, it seems worthwhile to note that the bound is rather subtle. An upper bound $\ll (\log P)^q$ would follow immediately by using H\"{o}lder's inequality to compare with the $q=1$ case, and applying standard Euler product calculations. Obtaining the crucial saving $1 + (1-q)\sqrt{\log\log P}$ in the denominator requires a careful analysis of the behaviour of various subproducts of $F_{P}^{\text{rand}}(1/2 + it)$. It is also shown in section 5.2 of Harper~\cite{harperrmflowmoments} that the upper bound in Multiplicative Chaos Result \ref{mcres1} is sharp, but we won't need (or be able) to exploit that here.

\vspace{12pt}
When dealing with character sums we cannot perform such a precise and complete ``conditioning'' procedure as in the genuine random multiplicative case, so it will be more straightforward (although not absolutely essential) to work with discrete sums rather than integral averages $\int_{-1/2}^{1/2}$ in the argument. We now deduce a discrete version of Multiplicative Chaos Result \ref{mcres1} to fit with this proof structure.

\begin{lem2}\label{disccontlem}
For any large $P$, we have
$$ \E \sum_{|k| \leq \frac{\log^{1.01}P}{2}} \int_{-\frac{1}{2\log^{1.01}P}}^{\frac{1}{2\log^{1.01}P}} |F_{P}^{\text{rand}}(1/2 + i\frac{k}{\log^{1.01}P} + it) - F_{P}^{\text{rand}}(1/2 + i\frac{k}{\log^{1.01}P})|^2 dt \ll \log^{0.99}P . $$
\end{lem2}

\begin{proof}[Proof of Lemma \ref{disccontlem}]
Since we have $F_{P}^{\text{rand}}(s) = \sum_{\substack{n=1, \\  n \; \text{is} \; P \; \text{smooth}}}^{\infty} \frac{f(n)}{n^s}$, with $f(n)$ a Steinhaus random multiplicative function, a mean square calculation shows that the left hand side in Lemma \ref{disccontlem} is
\begin{eqnarray}
& = & \sum_{|k| \leq \frac{\log^{1.01}P}{2}} \int_{-\frac{1}{2\log^{1.01}P}}^{\frac{1}{2\log^{1.01}P}} \E |\sum_{\substack{n=1, \\  n \; \text{is} \; P \; \text{smooth}}}^{\infty} \frac{f(n)}{n^{1/2 + i\frac{k}{\log^{1.01}P}}} (n^{-it} - 1)|^2 dt \nonumber \\
& = & \sum_{|k| \leq \frac{\log^{1.01}P}{2}} \int_{-\frac{1}{2\log^{1.01}P}}^{\frac{1}{2\log^{1.01}P}} \sum_{\substack{n=1, \\  n \; \text{is} \; P \; \text{smooth}}}^{\infty} \frac{|n^{-it}-1|^2}{n} dt . \nonumber
\end{eqnarray}
Here we have $|n^{-it} - 1| \ll \min\{|t|\log n, 1\} \leq \min\{\frac{\log n}{\log^{1.01}P}, 1\}$. Thus the contribution to the series from those $n \leq P^{\log\log P}$ is $\ll \frac{(\log\log P)^2}{\log^{0.02}P} \sum_{\substack{n=1, \\  n \; \text{is} \; P \; \text{smooth}}}^{\infty} \frac{1}{n}$, and since $\sum_{\substack{n=1, \\  n \; \text{is} \; P \; \text{smooth}}}^{\infty} \frac{1}{n} = \prod_{p \leq P} (1 - \frac{1}{p})^{-1} \ll \log P$ this is all $\ll (\log\log P)^2 \log^{0.98}P$. The contribution to the series from those $n > P^{\log\log P}$ (where we use the trivial bound $|n^{-it} - 1| \ll 1$) is also acceptably small, namely
$$ \ll e^{-\log\log P} \sum_{\substack{n=1, \\  n \; \text{is} \; P \; \text{smooth}}}^{\infty} \frac{1}{n^{1-1/\log P}} = e^{-\log\log P} \prod_{p \leq P} (1 - \frac{1}{p^{1-1/\log P}})^{-1} \ll e^{-\log\log P} \log P = 1 . $$
\end{proof}

\begin{multchaos2}\label{mcres2}
Uniformly for all large $P$ and $2/3 \leq q \leq 1$, we have
$$  \E( \frac{1}{\log^{1.01}P} \sum_{|k| \leq (\log^{1.01}P)/2} |F_{P}^{\text{rand}}(1/2 + i\frac{k}{\log^{1.01}P})|^2 )^q \ll \Biggl(\frac{\log P}{1 + (1-q)\sqrt{\log\log P}}\Biggr)^{q} . $$
\end{multchaos2}

\begin{proof}[Proof of Multiplicative Chaos Result \ref{mcres2}]
To deduce this from Multiplicative Chaos Result \ref{mcres1}, it will suffice to prove a suitable upper bound for
$$ \E( \sum_{|k| \leq (\log^{1.01}P)/2} \int_{-1/(2\log^{1.01}P)}^{1/(2\log^{1.01}P)} |F_{P}^{\text{rand}}(1/2 + i\frac{k}{\log^{1.01}P} + it) - F_{P}^{\text{rand}}(1/2 + i\frac{k}{\log^{1.01}P})|^2 dt )^q . $$
But using H\"{o}lder's inequality to compare the $q$-th moment with the first moment, we see this quantity is at most
$$ \Biggl( \E \sum_{|k| \leq \frac{\log^{1.01}P}{2}} \int_{-\frac{1}{2\log^{1.01}P}}^{\frac{1}{2\log^{1.01}P}} |F_{P}^{\text{rand}}(1/2 + i\frac{k}{\log^{1.01}P} + it) - F_{P}^{\text{rand}}(1/2 + i\frac{k}{\log^{1.01}P})|^2 dt \Biggr)^q , $$
which we can bound acceptably using Lemma \ref{disccontlem}. 
\end{proof}

\section{Proof of Theorem \ref{mainthmchar}}\label{secproof1}

\subsection{Notation and set-up}\label{thm1small}
We may restrict attention to the range $2/3 \leq q \leq 1$, since if $q$ is smaller we can use H\"{o}lder's inequality to upper bound $\Echar |\sum_{n \leq x} \chi(n)|^{2q}$ by $(\Echar |\sum_{n \leq x} \chi(n)|^{4/3})^{3q/2}$, and invoke the $q=2/3$ case.

Recall that $L := \min\{x,r/x\}$, which we may assume to be large since otherwise the Theorem is trivial. We have a parameter $P$ at our disposal, which we must choose to be comparable to $L$ on a {\em doubly logarithmic} scale, but (it turns out) somewhat smaller than $L$ on a logarithmic scale. In fact it will suffice if $P$ is around $\exp\{\log^{1/6}L\}$, and (for very minor technical reasons) we shall actually choose $P$ to be the largest number below $\exp\{\log^{1/6}L\}$ such that $\log^{0.01}P$ is an integer. Thus we will have $\log P \asymp \log^{1/6}L$ and $\log\log P \asymp \log\log L$, and to prove Theorem \ref{mainthmchar} it will suffice to show that
\begin{equation}\label{stpmaindisplay}
\Echar |\sum_{n \leq x} \chi(n)|^{2q} \ll \Biggl(\frac{x}{1 + (1-q)\sqrt{\log\log P}} \Biggr)^q .
\end{equation}

We set $M := 2\log^{1.02}P$, say (note this is an integer with our choice of $P$), and for each integer $k$ satisfying $|k| \leq M$ and each character $\chi$ mod $r$ we set $S_{k}(\chi) := \Re \sum_{p \leq P} (\frac{\chi(p)}{p^{1/2 + ik/\log^{1.01}P}} + \frac{\chi(p)^2}{2p^{1 + 2ik/\log^{1.01}P}})$. These are the prime number sums on which we shall ``condition'' in the next subsection.

Finally, recall that in Approximation Result \ref{apres} we have two parameters $N, \delta$ to set, for our construction of functions $g_{j}(\cdot)$ forming a smooth partition of unity. We shall make the final choices of these at the end of the proof, but it will turn out that the ``precision'' parameter $\delta$ may be chosen as a suitable negative power of $\log P$, and the ``range'' parameter $N$ (for which we have much flexibility) as a suitable multiple of $\log\log P$, say.

\subsection{The conditioning argument}\label{subsecmaincond}
Firstly, let $P(n)$ denote the largest prime factor of $n$, and as usual set $\Psi(x,y) := \#\{n \leq x : n \; \text{is} \; y \; \text{smooth}\} = \#\{n \leq x : P(n) \leq y \}$. Then we have
$$ \Echar |\sum_{n \leq x, P(n) \leq x^{1/\log\log x}} \chi(n)|^{2q} \leq \Biggl(\Echar |\sum_{n \leq x, P(n) \leq x^{1/\log\log x}} \chi(n)|^{2} \Biggr)^q = \Psi(x,x^{1/\log\log x})^q , $$
by H\"{o}lder's inequality and orthogonality of Dirichlet characters. Using standard smooth number estimates (see Theorem 7.6 of Montgomery and Vaughan~\cite{mv}, for example) this is $\ll (x (\log x)^{-c\log\log\log x})^q$, which is a negligible contribution in Theorem \ref{mainthmchar}.

So it will suffice to bound $\Echar |\sum_{n \leq x, P(n) > x^{1/\log\log x}} \chi(n)|^{2q}$. In the random version of the argument, one proceeds by conditioning on the values of $f(p)$ on all ``small'' primes $p$. We now look to emulate this procedure for character sums, breaking up the average $\Echar$ according to the values of all of the sums $S_{k}(\chi)$. Thus, recalling that the $g_j$ from Approximation Result \ref{apres} form a partition of unity, we may rewrite $\Echar |\sum_{n \leq x, P(n) > x^{1/\log\log x}} \chi(n)|^{2q}$ as
\begin{eqnarray}
&& \Echar \prod_{i=-M}^{M} (\sum_{j=-N}^{N+1} g_{j}(S_{i}(\chi))) \Biggl|\sum_{\substack{n \leq x, \\ P(n) > x^{1/\log\log x}}} \chi(n) \Biggr|^{2q} \nonumber \\
& = & \sum_{-N \leq j(-M), ..., j(0), ..., j(M) \leq N+1} \Echar \prod_{i=-M}^{M} g_{j(i)}(S_{i}(\chi)) \Biggl|\sum_{\substack{n \leq x, \\ P(n) > x^{1/\log\log x}}} \chi(n) \Biggr|^{2q} \nonumber \\
& = & \sum_{-N \leq j(-M) , ... , j(0), ..., j(M) \leq N+1} \sigma(\textbf{j}) \E^{\textbf{j}} \Biggl|\sum_{\substack{n \leq x, \\ P(n) > x^{1/\log\log x}}} \chi(n) \Biggr|^{2q} , \nonumber
\end{eqnarray}
where for any $(2M+1)$-vector $\textbf{j}$ from the outer sum we set $\sigma(\textbf{j}) := \Echar \prod_{i=-M}^{M} g_{j(i)}(S_{i}(\chi))$, and $\E^{\textbf{j}} W := \sigma(\textbf{j})^{-1} \Echar W \prod_{i=-M}^{M} g_{j(i)}(S_{i}(\chi))$ for all functions $W(\chi)$. This $\E^{\textbf{j}}$ is our character sum version of a conditional expectation. In particular, for the constant function 1, by the definitions we have $\E^{\textbf{j}} 1 = 1$ for all choices of the vector $\textbf{j}$. So applying H\"older's inequality to $\E^{\textbf{j}}$, we conclude overall that
$$ \Echar \Biggl|\sum_{\substack{n \leq x, \\ P(n) > x^{1/\log\log x}}} \chi(n) \Biggr|^{2q} \leq \sum_{-N \leq j(-M) , ..., j(0), ..., j(M) \leq N+1} \sigma(\textbf{j})\Biggl( \E^{\textbf{j}} \Biggl|\sum_{\substack{n \leq x, \\ P(n) > x^{1/\log\log x}}} \chi(n) \Biggr|^{2} \Biggr)^{q} . $$

\subsection{Passing to the random case}\label{subsecmaingorandom}
At this point, we can use Propositions \ref{condexpprop} and \ref{boxprobprop} to move from working with $\sigma(\textbf{j})$ and $\E^{\textbf{j}} \Biggl|\sum_{\substack{n \leq x, \\ P(n) > x^{1/\log\log x}}} \chi(n) \Biggr|^{2}$ to working with their analogues involving random multiplicative functions. Actually it isn't essential to do this at such an early stage in the argument, but doing it early will simplify various steps of the analysis, and will ultimately allow us to establish \eqref{stpmaindisplay} once we reach a point where we can invoke our Multiplicative Chaos Results.

Using Proposition \ref{condexpprop} with $Y=2M+1$, {\em provided that our choices of $N, \delta$ ultimately satisfy $x P^{400((2M+1)/\delta)^2 \log(N\log P)} < r$} we will have
$$ \E^{\textbf{j}} \Biggl|\sum_{\substack{n \leq x, \\ P(n) > x^{1/\log\log x}}} \chi(n) \Biggr|^{2} = \frac{1}{\sigma(\textbf{j})} \Biggl( \E \prod_{i=-M}^{M} g_{j(i)}(S_{i}(f)) \Biggl|\sum_{\substack{n \leq x, \\ P(n) > x^{1/\log\log x}}} f(n) \Biggr|^{2} + O\left(\frac{x}{(N \log P)^{(2M+1)(1/\delta)^2}} \right) \Biggr) . $$
Now observe that $\sum_{\textbf{j}} \sigma(\textbf{j}) = \Echar \prod_{i=-M}^{M} (\sum_{j=-N}^{N+1} g_{j}(S_{i}(\chi))) = \Echar 1 = 1$. Thus, applying H\"{o}lder's inequality to $\sum_{\textbf{j}}$, we see the total contribution to $\Echar |\sum_{n \leq x, P(n) > x^{1/\log\log x}} \chi(n)|^{2q}$ from all of the ``big Oh'' terms from Proposition \ref{condexpprop} is
$$ \ll \Biggl( \sum_{-N \leq j(-M) , ..., j(0), ..., j(M) \leq N+1} \sigma(\textbf{j}) \cdot \frac{1}{\sigma(\textbf{j})} \frac{x}{(N \log P)^{(2M+1)(1/\delta)^2}} \Biggr)^{q} = \Biggl( x (\frac{(2N+2)}{(N \log P)^{(1/\delta)^2}})^{2M+1} \Biggr)^{q} . $$
This is negligible for \eqref{stpmaindisplay}. 

Next, if we define $\sigma^{\text{rand}}(\textbf{j}) := \E \prod_{i=-M}^{M} g_{j(i)}(S_{i}(f))$ for all $(2M+1)$-vectors $\textbf{j}$, where $f$ is a Steinhaus random multiplicative function, then using Proposition \ref{boxprobprop} we get $\sigma(\textbf{j})^{1-q} \ll \sigma^{\text{rand}}(\textbf{j})^{1-q} + (\frac{1}{(N \log P)^{(2M+1)(1/\delta)^2}})^{1-q}$. Therefore we have
\begin{eqnarray}
&& \sum_{-N \leq j(-M) , ..., j(0), ..., j(M) \leq N+1} \sigma(\textbf{j})\Biggl( \frac{1}{\sigma(\textbf{j})} \E \prod_{i=-M}^{M} g_{j(i)}(S_{i}(f)) \Biggl|\sum_{\substack{n \leq x, \\ P(n) > x^{1/\log\log x}}} f(n) \Biggr|^{2} \Biggr)^q \nonumber \\
& \ll & \sum_{\textbf{j}} \sigma^{\text{rand}}(\textbf{j})\Biggl( \frac{1}{\sigma^{\text{rand}}(\textbf{j})} \E \prod_{i=-M}^{M} g_{j(i)}(S_{i}(f)) \Biggl|\sum_{\substack{n \leq x, \\ P(n) > x^{1/\log\log x}}} f(n) \Biggr|^{2} \Biggr)^q + \nonumber \\
&& + (\frac{1}{(N \log P)^{(2M+1)(1/\delta)^2}})^{1-q} \sum_{\textbf{j}} \Biggl( \E \prod_{i=-M}^{M} g_{j(i)}(S_{i}(f)) \Biggl|\sum_{\substack{n \leq x, \\ P(n) > x^{1/\log\log x}}} f(n) \Biggr|^{2} \Biggr)^q . \nonumber
\end{eqnarray}
Yet another application of H\"older's inequality to the sum over $\textbf{j}$, and recalling again that the $g_j$ form a partition of unity, implies that the final line here is
\begin{eqnarray}
& \ll & (\frac{1}{(N \log P)^{(2M+1)(1/\delta)^2}})^{1-q} \cdot ((2N+2)^{2M+1})^{1-q} \cdot \Biggl( \sum_{\textbf{j}} \E \prod_{i=-M}^{M} g_{j(i)}(S_{i}(f)) \Biggl|\sum_{\substack{n \leq x, \\ P(n) > x^{1/\log\log x}}} f(n) \Biggr|^{2} \Biggr)^q \nonumber \\
& = & \Biggl( (\frac{(2N+2)}{(N \log P)^{(1/\delta)^2}})^{2M+1} \Biggr)^{1-q} \Biggl( \E \Biggl|\sum_{\substack{n \leq x, \\ P(n) > x^{1/\log\log x}}} f(n) \Biggr|^{2} \Biggr)^q . \nonumber
\end{eqnarray}
This is certainly $\leq e^{-(1-q)\log\log P} x^q$, say, which is acceptable for \eqref{stpmaindisplay}.

In summary, if we now define $\E^{\textbf{j}, \text{rand}} W := \sigma^{\text{rand}}(\textbf{j})^{-1} \E W \prod_{i=-M}^{M} g_{j(i)}(S_{i}(f))$ for all random variables $W$, then to prove Theorem \ref{mainthmchar} it remains for us to show that
\begin{equation}\label{stpnextdisplay}
\sum_{-N \leq j(-M) , ..., j(0), ..., j(M) \leq N+1} \sigma^{\text{rand}}(\textbf{j}) \Biggl( \E^{\textbf{j}, \text{rand}} \Biggl|\sum_{\substack{n \leq x, \\ P(n) > x^{1/\log\log x}}} f(n) \Biggr|^{2} \Biggr)^{q} \ll \Biggl(\frac{x}{1 + (1-q)\sqrt{\log\log P}} \Biggr)^q .
\end{equation}
Note that we have now removed all mention of Dirichlet characters, and \eqref{stpnextdisplay} is purely a statement about random multiplicative functions.

\subsection{Passing to Euler products}\label{subsecpasseuler}
Our next step is to move from the left hand side of \eqref{stpnextdisplay}, which still involves {\em sums} of $f(n)$, to an expression featuring Euler products. This part of the argument will be very similar to the corresponding work in section 2.4 of Harper~\cite{harperrmflowmoments}, although not exactly the same because the way we set up our ``conditioning'' (using a smooth partition of unity) was necessarily different here.

Firstly (and crucially), since the $f(p)$ are {\em independent} mean zero random variables, and the various expressions $\prod_{i=-M}^{M} g_{j(i)}(S_{i}(f))$ in the definition of $\E^{\textbf{j}, \text{rand}}$ only involve the $f(p)$ for $p \leq P$ (not those for $p > P$), we find by expanding the square that
$$ \sum_{\textbf{j}} \sigma^{\text{rand}}(\textbf{j}) \Biggl( \E^{\textbf{j}, \text{rand}} \Biggl|\sum_{\substack{n \leq x, \\ P(n) > x^{1/\log\log x}}} f(n) \Biggr|^{2} \Biggr)^{q} = \sum_{\textbf{j}} \sigma^{\text{rand}}(\textbf{j}) \Biggl( \E^{\textbf{j}, \text{rand}} \sum_{\substack{m \leq x, \\ P(m) > x^{1/\log\log x} , \\ p|m \Rightarrow p > P}} \Biggl|\sum_{\substack{n \leq x/m, \\ n \; \text{is} \; P \; \text{smooth}}} f(n) \Biggr|^{2} \Biggr)^{q} . $$
Here we implicitly used the fact that $P < x^{1/\log\log x}$.

Setting $X = e^{\sqrt{\log x}}$, and replacing the condition that $P(m) > x^{1/\log\log x}$ by the weaker condition that $m > x^{1/\log\log x}$ (for an upper bound), and introducing an integral to smooth out on a scale of $1/X$, we find the bracketed term is
\begin{eqnarray}\label{smoothrandsplit}
& \ll & \Biggl(\E^{\textbf{j}, \text{rand}} \sum_{\substack{x^{1/\log\log x} < m \leq x, \\ p|m \Rightarrow p > P}} \frac{X}{m} \int_{m}^{m(1+1/X)} |\sum_{\substack{n \leq x/t, \\ n \; \text{is} \; P \; \text{smooth}}} f(n)|^{2} dt \Biggr)^q \nonumber \\
&& + \Biggl(\E^{\textbf{j}, \text{rand}} \sum_{\substack{x^{1/\log\log x} < m \leq x, \\ p|m \Rightarrow p > P}} \frac{X}{m} \int_{m}^{m(1+1/X)} |\sum_{\substack{x/t < n \leq x/m, \\ n \; \text{is} \; P \; \text{smooth}}} f(n)|^{2} dt \Biggr)^q .
\end{eqnarray}

Now observe again that $\sum_{\textbf{j}} \sigma^{\text{rand}}(\textbf{j}) = \E \prod_{i=-M}^{M} (\sum_{j=-N}^{N+1} g_{j}(S_{i}(f))) = \E 1 = 1$. Thus, applying H\"{o}lder's inequality to the sum over $\textbf{j}$, the total contribution to \eqref{stpnextdisplay} from the second bracket in \eqref{smoothrandsplit} is
\begin{eqnarray}
& \ll & \sum_{\textbf{j}} \sigma^{\text{rand}}(\textbf{j}) \Biggl( \E^{\textbf{j}, \text{rand}} \sum_{\substack{x^{1/\log\log x} < m \leq x, \\ p|m \Rightarrow p > P}} \frac{X}{m} \int_{m}^{m(1+1/X)} |\sum_{\substack{x/t < n \leq x/m, \\ n \; \text{is} \; P \; \text{smooth}}} f(n)|^{2} dt \Biggr)^q \nonumber \\
& \leq & \Biggl(\sum_{\textbf{j}} \sigma^{\text{rand}}(\textbf{j}) \E^{\textbf{j}, \text{rand}} \sum_{\substack{x^{1/\log\log x} < m \leq x, \\ p|m \Rightarrow p > P}} \frac{X}{m} \int_{m}^{m(1+1/X)} |\sum_{\substack{x/t < n \leq x/m, \\ n \; \text{is} \; P \; \text{smooth}}} f(n)|^{2} dt \Biggr)^q . \nonumber
\end{eqnarray}
But recalling the definitions of everything, and then using the orthogonality of random multiplicative functions, we see this is the same as
$$ \Biggl(\E \sum_{\substack{x^{1/\log\log x} < m \leq x, \\ p|m \Rightarrow p > P}} \frac{X}{m} \int_{m}^{m(1+1/X)} |\sum_{\substack{x/t < n \leq x/m, \\ n \; \text{is} \; P \; \text{smooth}}} f(n)|^{2} dt \Biggr)^q \ll \Biggl(\sum_{\substack{x^{1/\log\log x} < m \leq x, \\ p|m \Rightarrow p > P}} (1 + \frac{x}{mX}) \Biggr)^q . $$
Applying a standard sieve bound (e.g. Theorem 3.6 of Montgomery and Vaughan~\cite{mv}) to detect the condition that $p|m \Rightarrow p > P$ implies this quantity is $\ll (\frac{x}{\log P})^q$, which is more than good enough for us.

Meanwhile, the total contribution to  \eqref{stpnextdisplay} from the first bracket in \eqref{smoothrandsplit} is
\begin{eqnarray}
& \ll & \sum_{\textbf{j}} \sigma^{\text{rand}}(\textbf{j}) \Biggl(\E^{\textbf{j}, \text{rand}} \int_{x^{1/\log\log x}}^{x} |\sum_{\substack{n \leq x/t, \\ n \; \text{is} \; P \; \text{smooth}}} f(n)|^2 \sum_{\substack{t/(1+1/X) < m \leq t, \\ p|m \Rightarrow p > P}} \frac{X}{m} dt \Biggr)^q \nonumber \\
& \ll & \sum_{\textbf{j}} \sigma^{\text{rand}}(\textbf{j}) (\frac{1}{\log P})^q \Biggl(\E^{\textbf{j}, \text{rand}} \int_{x^{1/\log\log x}}^{x} |\sum_{\substack{n \leq x/t, \\ n \; \text{is} \; P \; \text{smooth}}} f(n)|^2 dt \Biggr)^q \nonumber \\
& = & \sum_{\textbf{j}} \sigma^{\text{rand}}(\textbf{j}) (\frac{x}{\log P})^q \Biggl(\E^{\textbf{j}, \text{rand}} \int_{1}^{x^{1-1/\log\log x}} |\sum_{\substack{n \leq z, \\ n \; \text{is} \; P \; \text{smooth}}} f(n)|^2 \frac{dz}{z^2} \Biggr)^q . \nonumber
\end{eqnarray}
Here the second line follows by applying a standard sieve bound again to the sum over $m$, and the final equality follows from the substitution $z=x/t$. Using Harmonic Analysis Result \ref{harmdirichlet}, we deduce that this is all
\begin{equation}\label{produpperintrand}
\ll (\frac{x}{\log P})^q \sum_{\textbf{j}} \sigma^{\text{rand}}(\textbf{j}) (\E^{\textbf{j}, \text{rand}} \int_{-\infty}^{\infty} \frac{|F_{P}^{\text{rand}}(1/2 + it)|^2}{|1/2+it|^2} dt)^q ,
\end{equation}
where $F_{P}^{\text{rand}}(s) := \prod_{p \leq P} (1 - \frac{f(p)}{p^s})^{-1} = \sum_{\substack{n=1, \\  n \; \text{is} \; P \; \text{smooth}}}^{\infty} \frac{f(n)}{n^s}$ is the Euler product corresponding to the random multiplicative function $f(n)$ on $P$-smooth numbers.

Since we have restricted to the range $2/3 \leq q \leq 1$, we can break up the integral in \eqref{produpperintrand} into sub-intervals of length 1 and obtain a bound
\begin{eqnarray}
& \leq & (\frac{x}{\log P})^q \sum_{\textbf{j}} \sigma^{\text{rand}}(\textbf{j}) \sum_{v=-\infty}^{\infty} (\E^{\textbf{j}, \text{rand}} \int_{v-1/2}^{v+1/2} \frac{|F_{P}^{\text{rand}}(1/2 + it)|^2}{|1/2+it|^2} dt)^q \nonumber \\
& \ll & (\frac{x}{\log P})^q \sum_{v=-\infty}^{\infty} \frac{1}{(|v|+1)^{4/3}} \sum_{\textbf{j}} \sigma^{\text{rand}}(\textbf{j}) (\E^{\textbf{j}, \text{rand}} \int_{v-1/2}^{v+1/2} |F_{P}^{\text{rand}}(1/2 + it)|^2 dt)^q . \nonumber
\end{eqnarray}
Those terms with $|v| > \log^{0.01}P$, say, trivially make a negligible contribution here. Indeed, using H\"{o}lder's inequality again their contribution is
\begin{eqnarray}
& \leq & (\frac{x}{\log P})^q \sum_{|v| > \log^{0.01}P} \frac{1}{(|v|+1)^{4/3}} (\sum_{\textbf{j}} \sigma^{\text{rand}}(\textbf{j}) \E^{\textbf{j}, \text{rand}} \int_{v-1/2}^{v+1/2} |F_{P}^{\text{rand}}(1/2 + it)|^2 dt)^q \nonumber \\
& = & (\frac{x}{\log P})^q \sum_{|v| > \log^{0.01}P} \frac{1}{(|v|+1)^{4/3}} (\E \int_{v-1/2}^{v+1/2} |F_{P}^{\text{rand}}(1/2 + it)|^2 dt)^q , \nonumber
\end{eqnarray}
and the orthogonality of random multiplicative functions implies this is
$$ = (\frac{x}{\log P})^q \sum_{|v| > \log^{0.01}P} \frac{1}{(|v|+1)^{4/3}} ( \sum_{\substack{n = 1, \\ n \; \text{is} \; P \; \text{smooth}}}^{\infty} \frac{1}{n} )^q \ll (\frac{x}{\log P})^q \frac{1}{\log^{1/300}P} \log^{q}P , $$
which is acceptable for Theorem \ref{mainthmchar}.

To complete the proof of the theorem, in view of \eqref{stpnextdisplay} and \eqref{produpperintrand} it will now suffice to show that uniformly for all $|v| \leq \log^{0.01}P$, we have
\begin{equation}\label{tobeprovedeqrand}
\sum_{\textbf{j}} \sigma^{\text{rand}}(\textbf{j}) (\E^{\textbf{j}, \text{rand}} \int_{v-1/2}^{v+1/2} |F_{P}^{\text{rand}}(1/2 + it)|^2 dt)^q \ll \Biggl(\frac{\log P}{1 + (1-q)\sqrt{\log\log P}}\Biggr)^{q} .
\end{equation}

\subsection{Strengthening the conditioning}\label{subsecrefinecond}
We now embark on establishing \eqref{tobeprovedeqrand}. For notational simplicity, we will write out the details in the case where $v=0$. The treatment of all other $|v| \leq \log^{0.01}P$ is exactly similar\footnote{Observe that our ``conditioning'' in $\E^{\textbf{j}, \text{rand}}$ is on the $S_{k}(f)$ for all $|k| \leq M = 2\log^{1.02}P$, corresponding to imaginary parts $\leq \frac{M}{\log^{1.01}P} = 2\log^{0.01}P$, which includes (with a little room to spare) the full range $|v| \leq \log^{0.01}P$ that we must handle in \eqref{tobeprovedeqrand}. Notice also that since we assume that $\log^{0.01}P$ is an integer, the points $i(v + \frac{k'}{\log^{1.01}P})$ for $v,k' \in \Z$ are of the desired form $i \frac{k}{\log^{1.01}P}$ with $k \in \Z$, which appear inside $S_{k}(f)$.}.

Firstly, since our ``conditional expectations'' $\E^{\textbf{j}, \text{rand}}$ contain information about all the sums $S_{k}(f)$, corresponding to the discrete points $t = \frac{k}{\log^{1.01}P}$, we need to move from the integral in \eqref{tobeprovedeqrand} to a discretised version. Thus we have
\begin{eqnarray}
&& \sum_{\textbf{j}} \sigma^{\text{rand}}(\textbf{j}) (\E^{\textbf{j}, \text{rand}} \int_{-1/2}^{1/2} |F_{P}^{\text{rand}}(1/2 + it)|^2 dt)^q \nonumber \\
& \ll & \sum_{\textbf{j}} \sigma^{\text{rand}}(\textbf{j}) ( \E^{\textbf{j}, \text{rand}} \sum_{|k| \leq \frac{\log^{1.01}P}{2}} \int_{-\frac{1}{2\log^{1.01}P}}^{\frac{1}{2\log^{1.01}P}} |F_{P}^{\text{rand}}(\frac{1}{2} + i\frac{k}{\log^{1.01}P}+it) - F_{P}^{\text{rand}}(\frac{1}{2} + i\frac{k}{\log^{1.01}P})|^2 dt )^q \nonumber \\
&& + \sum_{\textbf{j}} \sigma^{\text{rand}}(\textbf{j}) (\E^{\textbf{j}, \text{rand}} \frac{1}{\log^{1.01}P} \sum_{|k| \leq (\log^{1.01}P)/2} |F_{P}^{\text{rand}}(1/2 + i\frac{k}{\log^{1.01}P})|^2 )^q . \nonumber
\end{eqnarray}
Applying H\"older's inequality to the sum over $\textbf{j}$ once again, we see the first line here is
$$ \leq \Biggl( \E \sum_{|k| \leq \frac{\log^{1.01}P}{2}} \int_{-\frac{1}{2\log^{1.01}P}}^{\frac{1}{2\log^{1.01}P}} |F_{P}^{\text{rand}}(1/2 + i\frac{k}{\log^{1.01}P} + it) - F_{P}^{\text{rand}}(1/2 + i\frac{k}{\log^{1.01}P})|^2 dt \Biggr)^q , $$
which is $\ll \log^{0.99q}P$ in view of Lemma \ref{disccontlem} from section \ref{subsecrandeuler}. This is more than acceptable for \eqref{tobeprovedeqrand}, so it remains to handle the discrete sum on the second line.

It will also shortly be helpful to have some size restrictions on the $|F_{P}^{\text{rand}}(1/2 + i\frac{k}{\log^{1.01}P})|$, to aid us with setting the ``range'' parameter $N$ from Approximation Result \ref{apres}. So let us define the (random) ``bad set''
$$ \mathcal{T} := \{k \in \Z : |F_{P}^{\text{rand}}(1/2 + i\frac{k}{\log^{1.01}P})| \geq \log^{1.1}P \;\;\; \text{or} \;\;\; |F_{P}^{\text{rand}}(1/2 + i\frac{k}{\log^{1.01}P})| \leq \frac{1}{\log^{1.1}P} \} , $$
say. Splitting up the sum over $|k| \leq (\log^{1.01}P)/2$ according to whether $k \in \mathcal{T}$ or not, and using H\"{o}lder's inequality as (many times) before, we find that
\begin{eqnarray}
&& \sum_{\textbf{j}} \sigma^{\text{rand}}(\textbf{j}) (\E^{\textbf{j}, \text{rand}} \frac{1}{\log^{1.01}P} \sum_{|k| \leq (\log^{1.01}P)/2} |F_{P}^{\text{rand}}(1/2 + i\frac{k}{\log^{1.01}P})|^2 )^q \nonumber \\
& \leq & \sum_{\textbf{j}} \sigma^{\text{rand}}(\textbf{j}) (\E^{\textbf{j}, \text{rand}} \frac{1}{\log^{1.01}P} \sum_{\substack{|k| \leq (\log^{1.01}P)/2, \\ k \notin \mathcal{T}}} |F_{P}^{\text{rand}}(1/2 + i\frac{k}{\log^{1.01}P})|^2 )^q + \nonumber \\
&& + \Biggl( \frac{1}{\log^{1.01}P} \sum_{|k| \leq (\log^{1.01}P)/2} \E \textbf{1}_{k \in \mathcal{T}} |F_{P}^{\text{rand}}(1/2 + i\frac{k}{\log^{1.01}P})|^2 \Biggr)^q . \nonumber
\end{eqnarray}
On the final line, $\E \textbf{1}_{k \in \mathcal{T}} |F_{P}^{\text{rand}}(1/2 + i\frac{k}{\log^{1.01}P})|^2$ is at most $(\log^{1.1}P)^{-0.2} \E |F_{P}^{\text{rand}}(1/2 + i\frac{k}{\log^{1.01}P})|^{2.2} + (\frac{1}{\log^{1.1}P})^2$ (say). Standard results on the moments of random Euler products (see e.g. Euler Product Result 1 of \cite{harperrmfhigh}) imply that $(\log^{1.1}P)^{-0.2} \E |F_{P}^{\text{rand}}(1/2 + i\frac{k}{\log^{1.01}P})|^{2.2} \ll (\log^{1.1}P)^{-0.2} \log^{1.21}P = \log^{0.99}P$, so we get another more than acceptable contribution $\ll \log^{0.99q}P$ to \eqref{tobeprovedeqrand}.

Now the purpose of introducing the functions $g_{j(i)}$ (in the definitions of $\E^{\textbf{j}, \text{rand}}$ and $\sigma^{\text{rand}}(\textbf{j})$) was to approximately localise the values of the sums $S_{i}(f)$. At this stage, we will replace this approximate localisation by an exact version, with a view to ultimately producing a connection with Multiplicative Chaos Result \ref{mcres2}. Thus the contribution to \eqref{tobeprovedeqrand} from the sum over $k \notin \mathcal{T}$ is
\begin{eqnarray}
& \leq & \sum_{\textbf{j}} \sigma^{\text{rand}}(\textbf{j}) (\E^{\textbf{j}, \text{rand}} \frac{1}{\log^{1.01}P} \sum_{\substack{|k| \leq (\log^{1.01}P)/2, \\ k \notin \mathcal{T}}} \textbf{1}_{|S_{k}(f) - j(k)| \leq 1} |F_{P}^{\text{rand}}(1/2 + i\frac{k}{\log^{1.01}P})|^2 )^q + \nonumber \\
&& + \sum_{\textbf{j}} \sigma^{\text{rand}}(\textbf{j}) (\E^{\textbf{j}, \text{rand}} \frac{1}{\log^{1.01}P} \sum_{|k| \leq (\log^{1.01}P)/2} \textbf{1}_{|S_{k}(f) - j(k)| > 1} \textbf{1}_{k \notin \mathcal{T}} |F_{P}^{\text{rand}}(1/2 + i\frac{k}{\log^{1.01}P})|^2 )^q . \nonumber
\end{eqnarray}
Applying H\"{o}lder's inequality yet again to the sum over $\textbf{j}$ on the second line, and recalling the definitions of $\E^{\textbf{j}, \text{rand}}$ and $\sigma^{\text{rand}}(\textbf{j})$, we can bound that line by
$$ (\frac{1}{\log^{1.01}P} \sum_{|k| \leq \frac{\log^{1.01}P}{2}} \sum_{\textbf{j}} \E \prod_{i=-M}^{M} g_{j(i)}(S_{i}(f)) \cdot \textbf{1}_{|S_{k}(f) - j(k)| > 1} \textbf{1}_{k \notin \mathcal{T}} |F_{P}^{\text{rand}}(1/2 + i\frac{k}{\log^{1.01}P})|^2 )^q . $$
Now when $-N \leq j(k) \leq N$, by property (ii) from Approximation Result \ref{apres} we have $g_{j(k)}(S_{k}(f)) \cdot \textbf{1}_{|S_{k}(f) - j(k)| > 1} \leq \delta$. Furthermore, if $k \notin \mathcal{T}$ then $|S_{k}(f)| = |\log|F_{P}^{\text{rand}}(1/2 + i\frac{k}{\log^{1.01}P})|| + O(1) \leq 1.1\log\log P + O(1)$. Then {\em provided we take $N \geq 1.2\log\log P$ (say)}, when $j(k) = N+1$ we have $g_{j(k)}(S_{k}(f)) \cdot \textbf{1}_{k \notin \mathcal{T}} \leq \delta$, by property (iii) from Approximation Result \ref{apres}. So in any case we have $\prod_{i=-M}^{M} g_{j(i)}(S_{i}(f)) \cdot \textbf{1}_{|S_{k}(f) - j(k)| > 1} \textbf{1}_{k \notin \mathcal{T}} \leq \delta \prod_{i \neq k} g_{j(i)}(S_{i}(f))$, and then if we perform the sum over all possible values of $j(i)$ for $i \neq k$ we get $\delta \prod_{i \neq k} \sum_{-N \leq j(i) \leq N+1} g_{j(i)}(S_{i}(f)) = \delta$. Thus the second line above is at most
\begin{eqnarray}
(\frac{\delta}{\log^{1.01}P} \sum_{|k| \leq \frac{\log^{1.01}P}{2}} \sum_{-N \leq j(k) \leq N+1} \E |F_{P}^{\text{rand}}(1/2 + i\frac{k}{\log^{1.01}P})|^2 )^q & \ll & (\delta N \sum_{\substack{n = 1, \\ n \; \text{is} \; P \; \text{smooth}}}^{\infty} \frac{1}{n} )^q \nonumber \\
& \ll & (\delta N \log P)^q . \nonumber
\end{eqnarray}
This bound will be acceptable for \eqref{tobeprovedeqrand} {\em provided that $\delta \leq \frac{1}{N\sqrt{\log\log P}}$}.

\subsection{Conclusion}\label{subsecconclusion}
Finally, to establish \eqref{tobeprovedeqrand} (in the case $v=0$) it remains to prove that
\begin{eqnarray}
&& \sum_{\textbf{j}} \sigma^{\text{rand}}(\textbf{j}) \Biggl(\E^{\textbf{j}, \text{rand}} \frac{1}{\log^{1.01}P} \sum_{\substack{|k| \leq (\log^{1.01}P)/2, \\ k \notin \mathcal{T}}} \textbf{1}_{|S_{k}(f) - j(k)| \leq 1} |F_{P}^{\text{rand}}(1/2 + i\frac{k}{\log^{1.01}P})|^2 \Biggr)^q \nonumber \\
& \ll & \Biggl(\frac{\log P}{1 + (1-q)\sqrt{\log\log P}}\Biggr)^{q} . \nonumber
\end{eqnarray}

Recall now that $|F_{P}^{\text{rand}}(1/2 + i\frac{k}{\log^{1.01}P})| = \exp\{-\Re \sum_{p \leq P} \log(1 - \frac{f(p)}{p^{1/2 + i\frac{k}{\log^{1.01}P}}})\} \asymp \exp\{S_{k}(f)\}$ for all $k$ and all realisations of the random multiplicative function $f(n)$. And the point of our manipulations has been that the only $k$ values that now contribute to the sum over $k$ are those for which $S_{k}(f) \in [j(k)-1,j(k)+1]$. Noting also that we must have $|S_{k}(f)| \leq 1.1\log\log P + O(1)$ if $k \notin \mathcal{T}$, we see the left hand side is
\begin{eqnarray}
& \ll & \sum_{\textbf{j}} \sigma^{\text{rand}}(\textbf{j}) (\E^{\textbf{j}, \text{rand}} \frac{1}{\log^{1.01}P} \sum_{\substack{|k| \leq (\log^{1.01}P)/2, \\ |j(k)| \leq 1.1\log\log P + O(1)}} e^{2j(k)} )^q \nonumber \\
& = & \sum_{\textbf{j}} \sigma^{\text{rand}}(\textbf{j}) (\frac{1}{\log^{1.01}P} \sum_{\substack{|k| \leq (\log^{1.01}P)/2, \\ |j(k)| \leq 1.1\log\log P + O(1)}} e^{2j(k)} )^q . \nonumber
\end{eqnarray}
Here the quantity inside the bracket is a deterministic function of $\textbf{j}$, with the only remaining appearance of any randomness coming inside the multipliers $\sigma^{\text{rand}}(\textbf{j})$ in the outer ``averaging'' over $\textbf{j}$. Thus we are very close to the structure of the quantity bounded in Multiplicative Chaos Result \ref{mcres2}, where all of the averaging $\E$ occurs on the outside of the $q$-th power.

Indeed, recalling the definition of $\sigma^{\text{rand}}(\textbf{j})$ we can rewrite the sum here as
\begin{eqnarray}
&& \E \sum_{\textbf{j}} \prod_{i=-M}^{M} g_{j(i)}(S_{i}(f)) (\frac{1}{\log^{1.01}P} \sum_{\substack{|k| \leq (\log^{1.01}P)/2, \\ |j(k)| \leq 1.1\log\log P + O(1)}} e^{2j(k)} )^q \nonumber \\
& \ll & \E \sum_{\textbf{j}} \prod_{i=-M}^{M} g_{j(i)}(S_{i}(f)) (\frac{1}{\log^{1.01}P} \sum_{\substack{|k| \leq (\log^{1.01}P)/2, \\ |j(k)| \leq 1.1\log\log P + O(1)}} |F_{P}^{\text{rand}}(1/2 + i\frac{k}{\log^{1.01}P})|^2 )^q + \nonumber \\
&& + \E \sum_{\textbf{j}} \prod_{i=-M}^{M} g_{j(i)}(S_{i}(f)) (\frac{1}{\log^{1.01}P} \sum_{\substack{|k| \leq (\log^{1.01}P)/2, \\ |j(k)| \leq 1.1\log\log P + O(1)}} \textbf{1}_{|S_{k}(f) - j(k)| > 1} \log^{2.2}P )^q . \nonumber
\end{eqnarray}
And since $\sum_{\textbf{j}} \prod_{i=-M}^{M} g_{j(i)}(S_{i}(f)) \equiv 1$, the first line on the right hand side is at most $\E( \frac{1}{\log^{1.01}P} \sum_{|k| \leq (\log^{1.01}P)/2} |F_{P}^{\text{rand}}(1/2 + i\frac{k}{\log^{1.01}P})|^2 )^q$, which is $\ll \left(\frac{\log P}{1 + (1-q)\sqrt{\log\log P}}\right)^{q}$ by Multiplicative Chaos Result \ref{mcres2}.

Meanwhile, one last application of H\"older's inequality (to the expectation and sum over $\textbf{j}$ simultaneously), and using properties (ii) and (iii) from Approximation Result \ref{apres} (as in section \ref{subsecrefinecond}) to deduce that $g_{j(k)}(S_{k}(f)) \textbf{1}_{|j(k)| \leq 1.1\log\log P + O(1)} \textbf{1}_{|S_{k}(f) - j(k)| > 1} \leq \delta$, reveals the second line is
\begin{eqnarray}
& \leq & \Biggl( \E \sum_{\textbf{j}} \prod_{i=-M}^{M} g_{j(i)}(S_{i}(f)) \frac{1}{\log^{1.01}P} \sum_{\substack{|k| \leq (\log^{1.01}P)/2, \\ |j(k)| \leq 1.1\log\log P + O(1)}} \textbf{1}_{|S_{k}(f) - j(k)| > 1} \log^{2.2}P \Biggr)^q \nonumber \\
& \leq & \left( \frac{1}{\log^{1.01}P} \sum_{|k| \leq (\log^{1.01}P)/2} \log^{2.2}P \cdot \E \sum_{\textbf{j}} \delta \prod_{i \neq k} g_{j(i)}(S_{i}(f)) \right)^q . \nonumber
\end{eqnarray}
Since the functions $g_j$ form a partition of unity, this is all
$$ \leq \left( \frac{1}{\log^{1.01}P} \sum_{|k| \leq (\log^{1.01}P)/2} \log^{2.2}P \cdot \delta \sum_{-N \leq j(k) \leq N+1} 1 \right)^q \ll (\delta N \log^{2.2}P)^q , $$
which will be acceptable for \eqref{tobeprovedeqrand} {\em provided that $\delta \leq \frac{1}{N \log^{1.2}P \sqrt{\log\log P}}$}.

To conclude, we need only check that we can make an acceptable choice of the parameters $N, \delta$. In section \ref{subsecmaingorandom}, we needed to have $x P^{400((2M+1)/\delta)^2 \log(N\log P)} < r$ so that Proposition \ref{condexpprop} could be applied. Recalling that $M = 2\log^{1.02}P$, we see this will be satisfied provided that $P^{6499 (\log^{2.04}P) (1/\delta)^2 \log(N\log P)} < r/x$, say. In section \ref{subsecrefinecond} we needed $N \geq 1.2\log\log P$ and $\delta \leq \frac{1}{N\sqrt{\log\log P}}$, and now we have the more stringent condition $\delta \leq \frac{1}{N \log^{1.2}P \sqrt{\log\log P}}$.

So taking $N = \lceil 1.2\log\log P \rceil$ and $\delta = \frac{1}{\log^{1.3}P}$, say, all of our conditions will be satisfied provided that $P^{6500 (\log^{4.64}P) \log\log P} < r/x$. This indeed holds with our choice of $P$, slightly smaller than $\exp\{\log^{1/6}L\}$.
\qed

\section{Proofs of the other Theorems and Corollaries}\label{secotherproofs}

\subsection{Proof of Theorem \ref{mainthmcont}}
The neatest way to proceed is to let $\Phi : \R \rightarrow \R$ be a fixed non-negative function that is $\geq \textbf{1}_{[0,1]}$, and whose Fourier transform $\widehat{\Phi}$ is supported in $[-1/2\pi,1/2\pi]$. For example, we can take $\Phi(x)$ to be a Beurling--Selberg function, as described in Vaaler's paper~\cite{vaaler}. Then if we write $\E^{\text{cont}} W(t)$ to denote the continuous average $\frac{1}{T \widehat{\Phi}(0)} \int_{-\infty}^{\infty} \Phi(t/T) W(t) dt$, we see
$$ \frac{1}{T} \int_{0}^{T} |\sum_{n \leq x} n^{it}|^{2q} dt \leq \widehat{\Phi}(0) \E^{\text{cont}} |\sum_{n \leq x} n^{it}|^{2q} \ll \E^{\text{cont}} |\sum_{n \leq x} n^{it}|^{2q} , $$
and so it will suffice to prove the claimed Theorem \ref{mainthmcont} bound for $\E^{\text{cont}} |\sum_{n \leq x} n^{it}|^{2q}$.

The point of introducing $\E^{\text{cont}}$ is that we have $\E^{\text{cont}} 1 = \frac{1}{T \widehat{\Phi}(0)} \int_{-\infty}^{\infty} \Phi(t/T) dt = \frac{1}{\widehat{\Phi}(0)} \int_{-\infty}^{\infty} \Phi(u) du = 1$, and crucially $\Econt U^{it} \overline{V^{it}} = \frac{1}{T \widehat{\Phi}(0)} \int_{-\infty}^{\infty} \Phi(t/T) e^{-it\log(V/U)} dt = \frac{\widehat{\Phi}((T/2\pi)\log(V/U))}{\widehat{\Phi}(0)} = \textbf{1}_{U=V}$ for all natural numbers $1 \leq U,V < T$. (If we worked with $\frac{1}{T} \int_{0}^{T}$ directly, there would be error terms rather than an exact equality here.) These are the same properties that we had for our character average $\Echar$ in the proof of Theorem \ref{mainthmchar}. In particular, the analogues of Lemma \ref{evenmomentlem} and of Propositions \ref{condexpprop} and \ref{boxprobprop} hold with $\Echar$ replaced by $\Econt$, with $\chi(n)$ replaced by $n^{it}$, and with the various conditions $xU^k , x P^{400(Y/\delta)^2 \log(N\log P)} , P^{400(Y/\delta)^2 \log(N\log P)} < r$ replaced by their obvious substitutes $xU^k , x P^{400(Y/\delta)^2 \log(N\log P)} , P^{400(Y/\delta)^2 \log(N\log P)} < T$. This means that we can repeat the argument in sections \ref{thm1small}--\ref{subsecmaingorandom}, simply replacing $L = \min\{x,r/x\}$ by $L_T := \min\{x,T/x\}$ and replacing $S_{k}(\chi)$ by $\Re \sum_{p \leq P} (\frac{p^{it}}{p^{1/2 + ik/\log^{1.01}P}} + \frac{p^{2it}}{2p^{1 + 2ik/\log^{1.01}P}})$. At the end of section \ref{subsecmaingorandom}, the proof of Theorem \ref{mainthmcont} then reduces to establishing exactly the same bound \eqref{stpnextdisplay} for random multiplicative functions that we already proved in sections \ref{subsecpasseuler}--\ref{subsecconclusion}.
\qed

\subsection{Proof of Corollary \ref{cordev}}
This is a simple consequence of Theorem \ref{mainthmchar} together with Markov's inequality.

Thus for any $0 \leq q \leq 1$, the left hand side in Corollary \ref{cordev} is
$$ \leq \frac{\Echar |\sum_{n \leq x} \chi(n)|^{2q}}{(\lambda\frac{\sqrt{x}}{(\log\log(10L))^{1/4}})^{2q}} . $$

If $\lambda \geq e^{\sqrt{\log\log(10L)}}$, then simply taking $q = 1$ and using the fact that $\Echar |\sum_{n \leq x} \chi(n)|^{2} \leq x$ yields the desired bound.

For smaller $\lambda$, if we set $q = 1 - \delta$ with $0 < \delta \leq 1$, and apply Theorem \ref{mainthmchar}, we obtain that
$$ \frac{\Echar |\sum_{n \leq x} \chi(n)|^{2q}}{(\lambda\frac{\sqrt{x}}{(\log\log(10L))^{1/4}})^{2q}} \ll \frac{1}{\lambda^2} \frac{\lambda^{2\delta} (\log\log(10L))^{q/2}}{(1 + \delta\sqrt{\log\log(10L)})^q} \leq \frac{1}{\lambda^2} \frac{\lambda^{2\delta}}{\delta} = \frac{\log\lambda}{\lambda^2} \frac{e^{2\delta \log\lambda}}{\delta \log\lambda}  . $$
Choosing $\delta = \frac{1}{2\log\lambda}$ yields the claimed upper bound.
\qed

\subsection{Proof of Corollary \ref{cortheta}}
As usual, it will suffice to prove Corollary \ref{cortheta} when $2/3 \leq q \leq 1$, because when $q$ is smaller we can use H\"{o}lder's inequality to compare with the $q=2/3$ case.

If $\chi$ is an even Dirichlet character mod $r$, then using the definition of the theta function and partial summation we have
\begin{eqnarray}
\theta(1,\chi) = \sum_{n=1}^{\infty} \chi(n) e^{-\pi n^{2}/r} & = & \sum_{n \leq \sqrt{r\log r}} \chi(n) e^{-\pi n^{2}/r} + O(1) \nonumber \\
& = & \sum_{n \leq \sqrt{r\log r}} \chi(n) \int_{n}^{\sqrt{r\log r}} \frac{2\pi u}{r} e^{-\pi u^{2}/r} du + O(1) \nonumber \\
& = & \int_{1}^{\sqrt{r\log r}} \frac{2\pi u}{r} e^{-\pi u^{2}/r} \sum_{n \leq u} \chi(n) du + O(1) . \nonumber
\end{eqnarray}
By noting that $\int_{1}^{\sqrt{r\log r}} \frac{2\pi u}{r} e^{-\pi u^{2}/r} du \asymp 1$, and applying H\"{o}lder's inequality to the integral over $u$ (here we use the fact that $2q \geq 1$), we deduce that for each even $\chi$ we have
$$ |\theta(1, \chi)|^{2q} \ll \int_{1}^{\sqrt{r\log r}} \frac{2\pi u}{r} e^{-\pi u^{2}/r} |\sum_{n \leq u} \chi(n)|^{2q} du + 1 . $$
The contribution to this integral from $1 \leq u \leq r^{1/4}$ is trivially $\ll \int_{1}^{r^{1/4}} \frac{u}{r} u^{2q} du \ll 1$, which is negligible. On the remaining range $r^{1/4} \leq u \leq \sqrt{r\log r}$, Theorem \ref{mainthmchar} implies that $\Echar |\sum_{n \leq u} \chi(n)|^{2q} \ll (\frac{u}{1+(1-q)\sqrt{\log\log r}})^q$, giving an acceptable contribution
$$ \ll (\frac{1}{1+(1-q)\sqrt{\log\log r}})^q \int_{r^{1/4}}^{\sqrt{r\log r}} \frac{u}{r} e^{-\pi u^{2}/r} u^{q} du \ll (\frac{\sqrt{r}}{1+(1-q)\sqrt{\log\log r}})^q $$
for Corollary \ref{cortheta}.

If $\chi$ is an odd character, then by definition we have $\theta(1,\chi) = \sum_{n=1}^{\infty} n\chi(n) e^{-\pi n^{2}/r}$, and a similar partial summation as before yields that
$$ \theta(1,\chi) = \int_{1}^{\sqrt{r\log r}} (\frac{2\pi u^2}{r} - 1) e^{-\pi u^{2}/r} \sum_{n \leq u} \chi(n) du + O(1) . $$
Since we now have $\int_{1}^{\sqrt{r\log r}} |\frac{2\pi u^2}{r} - 1| e^{-\pi u^{2}/r} du \asymp \sqrt{r}$, when we apply H\"older's inequality and Theorem \ref{mainthmchar} (as in the even character case) we now have an extra factor $(\sqrt{r})^{2q} = r^q$ in our bounds, producing the claimed upper bound for the average over odd characters.
\qed

\subsection{Proof of Theorem \ref{mainthmmulttwist}}
Once we introduce the multiplicative twist $h(n)$, the Euler product factor corresponding to a prime $p$ changes from being $(1 - \frac{\chi(p)}{p^{1/2+it}})^{-1}$ (or $(1 - \frac{f(p)}{p^{1/2+it}})^{-1}$ on the random multiplicative side), to
$$ 1 + \frac{h(p) \chi(p)}{p^{1/2+it}} + \sum_{k=2}^{\infty} \frac{h(p^k) \chi(p)^k}{p^{k(1/2+it)}} \;\;\;\;\; \text{or} \;\;\;\;\; 1 + \frac{h(p) f(p)}{p^{1/2+it}} + \sum_{k=2}^{\infty} \frac{h(p^k) f(p)^k}{p^{k(1/2+it)}} . $$
Here we have $|\frac{h(p) \chi(p)}{p^{1/2+it}} + \sum_{k=2}^{\infty} \frac{h(p^k) \chi(p)^k}{p^{k(1/2+it)}}| \leq \sum_{k=1}^{\infty} \frac{1}{p^{k/2}} = \frac{1}{\sqrt{p} - 1}$. Provided that $p \geq 5$, this is all $\leq \frac{1}{\sqrt{5} - 1} < 1$, so we can still apply Taylor expansion to analyse the logarithms of the Euler factors. This is not the case for the primes 2 and 3, but for those we still have an upper bound $\ll 1$ for the Euler factors.

With the above observations, the proof of Theorem \ref{mainthmmulttwist} is a fairly obvious adjustment of the proof of Theorem \ref{mainthmchar}. Rather than $S_{k}(\chi) := \Re \sum_{p \leq P} (\frac{\chi(p)}{p^{1/2 + ik/\log^{1.01}P}} + \frac{\chi(p)^2}{2p^{1 + 2ik/\log^{1.01}P}})$, in section \ref{secproof1} we work with $S_{k,h}(\chi) := \Re \sum_{5 \leq p \leq P} (\frac{h(p) \chi(p)}{p^{1/2 + ik/\log^{1.01}P}} + \frac{(h(p^2) - (1/2)h(p)^2)\chi(p)^2}{p^{1 + 2ik/\log^{1.01}P}})$ (coming from the first and second order terms in the Taylor expansions of the logarithms of the Euler factors). Note that the primes 2 and 3 are omitted here. The calculations then proceed essentially without change, until in \eqref{produpperintrand} we end up with $|\prod_{p \leq P} (1 + \frac{h(p) f(p)}{p^{1/2 +it}} + \sum_{k=2}^{\infty} \frac{h(p^k) f(p)^k}{p^{k(1/2 +it)}})|$ in place of $|F_{P}^{\text{rand}}(1/2+it)|$. And this is $\ll |F_{P, h}^{\text{rand}}(1/2+it)|$, where $F_{P, h}^{\text{rand}}(s) := \prod_{5 \leq p \leq P} (1 + \frac{h(p) f(p)}{p^{s}} + \sum_{k=2}^{\infty} \frac{h(p^k) f(p)^k}{p^{ks}})$ (with the primes 2 and 3 omitted).

Continuing through sections \ref{subsecrefinecond} and \ref{subsecconclusion}, we only need to verify that we have the same upper bound $\E |F_{P,h}^{\text{rand}}(1/2 + i\frac{k}{\log^{1.01}P})|^{2.2} \ll \log^{1.21}P$ as for $\E |F_{P}^{\text{rand}}(1/2 + i\frac{k}{\log^{1.01}P})|^{2.2}$, and (most importantly) that Multiplicative Chaos Result \ref{mcres1} continues to hold with $F_{P}^{\text{rand}}(1/2 + it)$ replaced by $F_{P,h}^{\text{rand}}(1/2 + it)$. The former is an easy modification of Euler Product Result 1 of \cite{harperrmfhigh}. The latter is also quite straightforward to check by working through the proofs in sections 3.1--3.2 and 4.1--4.3 of \cite{harperrmflowmoments}, noting that Lemma 1 from that paper holds without change for the Euler product $F_{P, h}^{\text{rand}}(s)$ (here it is important that $|h(p)|=1$ on primes $p$, which produce the main terms there), and all subsequent results flow from Lemma 1.
\qed


\end{document}